\documentclass[12pt]{article}
\usepackage{amsfonts,amssymb,amsmath,mathrsfs}
\usepackage{color}
\usepackage[colorlinks, linkcolor=blue, citecolor=blue]{hyperref}
\usepackage{graphicx}
\usepackage[T1]{fontenc}
\hypersetup{
    colorlinks=true,
    linkcolor=blue,
    filecolor=blue,
    urlcolor=blue,
    citecolor=blue
}

\allowdisplaybreaks[4]

\allowdisplaybreaks

\newtheorem{con0}{Theorem}[section]
\newtheorem{thm0}{Theorem}[section]
\newtheorem{exa0}{Theorem}[section]

\newtheorem{con1}[con0]{Condition}
\newtheorem{def1}[thm0]{Definition}
\newtheorem{lem1}[thm0]{Lemma}
\newtheorem{thm1}[thm0]{Theorem}
\newtheorem{cor1}[thm0]{Corollary}
\newtheorem{pro1}[thm0]{Proposition}
\newtheorem{rem1}[thm0]{Remark}
\newtheorem{ass1}[thm0]{Assumption}
\newtheorem{exa1}[exa0]{Example}

\def\blemma{\begin{lem1}}\def\elemma{\end{lem1}}
\def\btheorem{\begin{thm1}}\def\etheorem{\end{thm1}}

\def\bproposition{\begin{pro1}}\def\eproposition{\end{pro1}}

\def\benumerate{\begin{enumerate}}\def\eenumerate{\end{enumerate}}
\def\bitemize{\begin{itemize}}\def\eitemize{\end{itemize}}\def\itm{\item}

\def\beqlb{\begin{eqnarray}}\def\eeqlb{\end{eqnarray}}
\def\beqnn{\begin{eqnarray*}}\def\eeqnn{\end{eqnarray*}}
\def\beqa{\begin{equation}}\def\eeqa{\end{equation}}
\def\beqaa{\begin{equation*}}\def\eeqaa{\end{equation*}}

\def\bproof{\begin{proof}~}\def\eproof{\qed\end{proof}}
\def\qed{\hfill$\square$\smallskip}

\DeclareMathOperator{\var}{var}

\def\eqref#1{{\rm(\ref{#1})}}

\def\ar{\!\!\!&}\def\nnm{\nonumber}\def\ccr{\nnm\\}
\def\qqquad{\qquad\qquad}

\def\mrm{\mathrm}\def\mbf{\mathbf}\def\mcr{\mathscr}
\def\mbb{\mathbb}

\def\d{\mrm{d}}\def\e{\mrm{e}}\def\im{\mrm{i}}
\def\b{\mrm{b}}\def\var{\mrm{var}}

\def\itDelta{{\it\Delta}}

\def\itOmega{{\it\Omega}}

\def\const{\mrm{const}}

\def\blue{\color{blue}}

\begin{document}

\noindent{(Version: 2024-11-21)}

\bigskip\bigskip

\noindent{\bf\LARGE Derrida--Retaux type models and}

\smallskip

\noindent{\bf\LARGE related scaling limit theorems}

\bigskip

\noindent{Zenghu Li and Run Zhang\,\footnote{Corresponding author.}}

\medskip

\noindent{\it Laboratory of Mathematics and Complex Systems,}

\noindent{\it School of Mathematical Sciences,}

\noindent{\it Beijing Normal University, Beijing 100875, China}

\noindent{Emails: \tt lizh@bnu.edu.cn, zhangrun@mail.bnu.edu.cn}

\bigskip

\textbf{Abstract:~} We give characterizations of the transition semigroup and generator of a continuous-time Derrida--Retaux type process that generalizes the one introduced by Hu, Mallein and Pain (Commun. Math. Phys., 2020). It is shown that the process arises naturally as the scaling limit of the discrete-time max-type recursive models introduced by Hu and Shi (J. Stat. Phys., 2018).

\bigskip


\smallskip


\bigskip


\section{Introduction}\label{Sect1}

\setcounter{equation}{0}

A discrete-time max-type recursive model was introduced by Derrida and Retaux \cite{DeR14} in the study of the depinning transition in the limit of strong disorder. Write $T_a(x)= (x-a)_+$ for $a,x\ge 0$. For any function $f$ on $\mbb{R}_+:= [0,\infty)$, write
 \beqnn
T_af(x)= f\circ T_a (x)= f((x-a)_+), \quad a,x\ge 0.
 \eeqnn
Given a Borel measure $\mu$ on $\mbb{R}_+$, we denote by $\mu\circ T_a^{-1}$ the measure defined by
 \beqnn
\mu\circ T_a^{-1}(B)= \mu(\{x\ge 0: T_a(x)\in B\}), \quad B\in \mcr{B}(\mbb{R}_+),
 \eeqnn
where $\mcr{B}(\mbb{R_+})$ is the Borel $\sigma$-algebra on $\mbb{R}_+$. Then, given a probability measure $\mu_0$ on $\mbb{R}_+$, we can define a sequence of probability measures $(\mu_n: n\ge 0)$ recursively by
 \beqlb\label{DR0recur1mu}
\mu_{n+1}= (\mu_n^{*2})\circ T_1^{-1}, \quad n\ge 0,
 \eeqlb
where $\mu_n^{*2}= \mu_n*\mu_n$ denotes the convolution. The sequence $(\mu_n: n\ge 0)$ is called a \textit{Derrida--Retaux model}, or simply a \textit{DR model}. By \eqref{DR0recur1mu} it is easy to see that
 \beqnn
\int_{\mbb{R}_+} x \mu_{n+1}(\mathrm{d}x)
 =
\int_{\mbb{R}_+} (x-1)_+ \mu_n^{*2}(\mathrm{d}x)
 \le
2\int_{\mbb{R}_+} x \mu_n(\mathrm{d}x).
 \eeqnn
Therefore the decreasing limit exists:
 \beqnn
F_\infty= \lim_{n\to \infty} \,2^{-n}\!\int_{\mbb{R}_+} x \mu_n(\mathrm{d}x),
 \eeqnn
which is called the \textit{free energy}. The DR model is referred to as \textit{pinned} if $F_\infty> 0$, and as \textit{unpinned} if $F_\infty= 0$. One main problem in this study is to determine for which initial distribution $\mu_0$ the model is pinned or unpinned.

It is believed that for a large class of recursive models, including the DR model, there is a highly non-trivial phase transition. To discuss the phase transition from the pinned to the unpinned regime, it is convenient to specify the mass of $\mu_0$ at the origin. Consider the decomposition:
 \beqnn
\mu_0(\d x)= p\delta_0(\d x) + (1-p)\vartheta(\d x), \quad x\ge 0,
 \eeqnn
where $0\le p\le 1$ is a constant and $\vartheta$ is a fixed probability measure carried by $(0,\infty)$. Let $F_\infty(p)$ denote the associated free energy. Then $p\mapsto F_\infty(p)$ is a decreasing function on $[0,1]$. Write $p_c\in [0,1]$ for the \textit{critical parameter} distinguishing the pinned and the unpinned regimes, that is,
 \beqnn
p_c= \sup\{p\in [0,1]: F_\infty(p)> 0\}
 \eeqnn
with the convention $\sup\emptyset= 0$. Derrida and Retaux \cite{DeR14} conjectured that, under the assumption $p_c>0$ and some integrability conditions on $\vartheta$, there exists some constant $C>0$ such that
 \beqlb\label{DR0conj}
F_\infty(p)= \exp\Big(-\frac{C+o(1)}{\sqrt{p_c-p}}\Big), \quad p\uparrow p_c.
 \eeqlb
A weaker form of \eqref{DR0conj} has been proved by Chen et al.\ \cite{CDDHS21} in the special case where $\vartheta$ is carried by the set $\{1,2,\cdots\}$. Another basic question is the asymptotic behavior of the \textit{sustainability probability} $\mu_n(0,\infty)$ as $n\to \infty$. When $p= p_c$ and $\vartheta$ is carried by $\{1,2,\cdots\}$ it is expected that
 \beqlb\label{mu_n(0,infty)=4/n^2+o(.)}
\mu_n(\{1,2,\cdots\})= \frac{4}{n^2} + o\Big(\frac{1}{n^2}\Big), \quad n\to \infty.
 \eeqlb
We refer the reader to \cite{CDDHS21, CEGM84, DeR14} for the physical explanations of the above prediction.

A continuous-time version of the DR model was introduced by Hu et al.~\cite{HMP20}, who showed the model is exactly solvable and belongs to the universality class mentioned above. By definition, the model is a continuous-time flow of probability measures $(\mu_t: t\ge 0)$ on $\mbb{R}_+$ solving the differential equation:
 \beqlb\label{Cdifeq0a}
\partial_t\mu_t
 =
\mu_t^{*2} - \mu_t + \partial_x\mu_t 1_{\{x> 0\}}, \quad t\ge 0.
 \eeqlb
By \cite[Theorem~1.8]{HMP20}, for each initial state $\mu_0$ there is a unique weak solution to \eqref{Cdifeq0a}. Following Hu et al.~\cite{HMP20}, we call $(\mu_t: t\ge 0)$ a \textit{continuous-time DR model}, or simply a \textit{CDR model}. A differential equation similar to \eqref{Cdifeq0a} was informally derived by Derrida and Retaux \cite{DeR14} as the scaling limit of the model defined by \eqref{DR0recur1mu}, which has played the key role in the prediction \eqref{DR0conj}. The CDR model $(\mu_t: t\ge 0)$ is exactly solvable when it is started with the initial distribution
 \beqnn
\mu_0(\d x)= p\delta_0(\d x) + (1-p)\lambda\e^{-\lambda x}\d x, \quad x\ge 0,
 \eeqnn
where $0\le p\le 1$ and $\lambda>0$. In this case, the free energy is defined by
 \beqnn
F_\infty(p,\lambda)= \lim_{t\to \infty} \,\e^{-t}\!\int_{\mbb{R}_+} x \mu_t(\mathrm{d}x).
 \eeqnn
For the CDR model, Hu et al.~\cite{HMP20} characterized its pinned and unpinned classes of the parameters $(\lambda,p)$ and proved the Derrida--Retaux conjecture.

A discrete-time generalization of the DR model was introduced and studied by Hu and Shi \cite{HuS18}. Let $0\le \alpha\le 1$ and let $q= \{q_1,q_2,\cdots\}$ be a fixed discrete probability distribution on $\{1,2,\cdots\}$. Given a Borel probability measure $\mu$ on $\mbb{R}_+$, we define the measure $\mu^q$ by
 \beqlb\label{mu^q=def}
\mu^q= \sum_{k=1}^\infty q_k\mu^{*k},
 \eeqlb
where $\mu^{*k}$ denotes the $k$-fold convolution. The max-type model of Hu and Shi \cite{HuS18} can be defined by the recursive formula
 \beqlb\label{gDR0recur1mu}
\mu_{n+1}= [(1-\alpha)\mu_n + \alpha\mu_n*\mu_n^q]\circ T_1^{-1}, \quad n\ge 0.
 \eeqlb
It is natural to call $(\mu_n: n\ge 0)$ a \textit{generalized DR model} with \textit{renewal rate} $\alpha$ and \textit{offspring distribution} $q= \{q_1,q_2,\cdots\}$. When $\alpha= q_1= 1$, it reduces to the classical model \eqref{DR0recur1mu}. For the generalized DR model, Hu and Shi \cite{HuS18} showed a wide range for the exponent of the free energy in the nearly supercritical regime and Chen et al.~\cite{CHS22} established a weaker form of the conjecture \eqref{mu_n(0,infty)=4/n^2+o(.)}. A stronger result for the generalized DR-model with exponential-type marginal distributions was given by Li and Zhang \cite{LiZ24}.

In this work, we are interested in the scaling limits of the generalized DR model leading to continuities-time models like the one defined by \eqref{Cdifeq0a}. Let $(\mu_n: n\ge 0)$ be given by \eqref{gDR0recur1mu}. For $k\ge 1$ consider the rescaled measure $\gamma_n^{(k)}(\d x)= \mu_n(k\d x)$. From \eqref{gDR0recur1mu} it follows that
 \beqnn
\gamma_{n+1}^{(k)} - \gamma_n^{(k)}= \alpha[\gamma_n^{(k)}*(\gamma_n^{(k)})^q - \gamma_n^{(k)}]\circ T_{1/k}^{-1} + (\gamma_n^{(k)}\circ T_{1/k}^{-1} - \gamma_n^{(k)}).
 \eeqnn
Then, taking $\alpha= a/k$ for some $a\ge 0$, one naturally expects that the rescaled dynamics $(\gamma_{\lfloor kt\rfloor}^{(k)}: t\ge 0)$ would converge as $k\to \infty$ to the solution of
 \beqlb\label{gCDR0difeq0a}
\partial_t\mu_t
 =
a(\mu_t*\mu_t^q-\mu_t) + \partial_x\mu_t 1_{\{x> 0\}}, \quad t\ge 0,
 \eeqlb
where $\mu_t^q$ is defined as in \eqref{mu^q=def}; see \cite[p.280]{DeR14} and \cite[p.611]{HMP20}. This observation is confirmed by Theorem~\ref{thMDR-lim2a} of this paper, where the convergence of the probabilities in a Wasserstein distance is proved. We call the solution $(\mu_t: t\ge 0)$ of \eqref{gCDR0difeq0a} a \textit{generalized CDR model}.

Let $\b\mcr{C}(\mbb{R}_+)$ be the set of bounded continuous functions on $\mbb{R}_+$ and let $\b\mcr{C}^1(\mbb{R}_+)$ be the set of functions in $\b\mcr{C}(\mbb{R}_+)$ with bounded continuous first derivatives. For $t\ge 0$ and $f\in \b\mcr{C}^1(\mbb{R}_+)$ let
 \beqlb\label{CDRA_sf(x)=def}
A_tf(x)= a\int_{\mbb{R}_+}[f(x+z)-f(x)]\mu_t^q(\d z) - f'(x)1_{\{x> 0\}}, \quad x\ge 0.
 \eeqlb
Then the family of operators $(A_t: t\ge 0)$ generates an inhomogeneous transition semigroup $(P_{r,t}: t\ge r\ge 0)$ on $\mbb{R}_+$. We shall see that $(\mu_t: t\ge 0)$ is a \textit{closed entrance law} for $(P_{r,t}: t\ge r\ge 0)$, that is,
 \beqnn
\mu_t = \int_{\mbb{R}_+}\mu_0(\d x)P_{0,t}(x,\cdot), \quad t\ge 0.
 \eeqnn
If a positive Markov process $(X_t: t\ge 0)$ has transition semigroup $(P_{r,t}: t\ge r\ge 0)$, we call it a \textit{generalized CDR process} associated with the generalized CDR model $(\mu_t: t\ge 0)$. We shall see that $(X_t: t\ge 0)$ is a generalized CDR process if and only if, for every $f\in \b\mcr{C}^1(\mbb{R}_+)$,
 \beqlb\label{CDR-martprob1b}
f(X_t)= f(X_0) + \int_0^t A_sf(X_s)\d s + M_t(f), \quad t\ge 0,
 \eeqlb
where $\{M_t(f): t\ge 0\}$ is a martingale. In this case, if $X_0$ has distribution $\mu_0$, then $X_t$ has distribution $\mu_t$ for every $t\ge 0$.

Let $N(\d s,\d u)$ be a time-space Poisson random measure on $(0,\infty)\times (0,1)$ with intensity $a\d s\d u$. A c\'ad\'ag realization of the generalized CDR process is given by the pathwise unique solution to the stochastic integral equation:
 \beqlb\label{gCDRstoeq1X_t}
X_t= X_0 + \int_{(0,t]}\int_{(0,1)} G_s^{-1}(u) N(\d s,\d u) - \int_0^t 1_{\{X_s> 0\}} \d s, \quad t\ge 0,
 \eeqlb
where $G_s^{-1}$ denotes the right-continuous inverse of the distribution function of $\mu_s^q$. A special form of \eqref{gCDRstoeq1X_t} has been used by Hu et al.~\cite{HMP20} in their construction of the CDR process associated with the model defined by \eqref{Cdifeq0a}.

Suppose that $(\mu_n: n\ge 0)$ is a generalized DR model defined by \eqref{gDR0recur1mu}. Let $U_n, \eta_n,$ $n\ge 0$ be independent random variables, where the $U_n$ follows the uniform distribution $U(0,1)$ and the $\eta_n$ the Bernoulli distribution $B(1,\alpha)$, that is, $\mbf{P}(\eta_n=1)= \alpha$ and $\mbf{P}(\eta_n=0)= 1-\alpha$. Given a positive random variable $X_0$ independent of $(U_n,\eta_n: n\ge 0)$, define recursively
 \beqlb\label{MDR-recurs2X_n}
X_{n+1}= (X_n + \eta_nG_n^{-1}(U_n)-1)_+, \quad n\ge 0,
 \eeqlb
where $G_n^{-1}$ is the right-continuous inverse of the distribution function of $\mu_n^q$. We show in Theorem~\ref{thMDR-conv0Skor} of this paper that the generalized CDR process defined by \eqref{gCDRstoeq1X_t} arises naturally as the limit in the Skorokhod space of the rescaled sequence $(k^{-1}X_{\lfloor kt\rfloor}: t\ge 0)$ as $k\to \infty$.

The remainder of the paper is organized as follows. The basic properties of the generalized DR models with continuous-times and discrete-times are discussed in Section~2, where the limit theorem for the rescaled dynamics $(\gamma_{\lfloor kt\rfloor}^{(k)}: t\ge 0)$ is proved. In Section~3, we give characterizations of the transition semigroup and generator of the generalized CDR process. The martingale problem of the process is discussed in Section~4. The convergence of the rescaled process in the Skorokhod space is proved in Section~5.

\section{Derrida--Retaux type models}\label{Sect2}

\setcounter{equation}{0}

\subsection{Preliminaries}\label{SSct2.1}

Let $\b\mcr{B}(\mbb{R_+})$ be the set of bounded Borel functions on $\mbb{R}_+$. For $f\in \b\mcr{B}(\mbb{R_+})$, we define its \textit{supremum norm} $\|f\|_\infty= \sup_{x\in \mbb{R}_+} |f(x)|$ and its \textit{$\rho$-Lipschitz seminorm}
 \beqnn
\|f\|_\rho=\sup_{x\neq y\in \mbb{R}_+} \rho(x,y)^{-1}|f(x)-f(y)|,
 \eeqnn
where $\rho(x,y)= 1\land |x-y|$ denotes the \textit{truncated Euclidean distance}.

Let $\mcr{P}(\mbb{R}_+)$ be the space of Borel probability measures on $\mbb{R}_+$. For any $\mu,\nu\in \mcr{P}(\mbb{R}_+)$ let $\mcr{C}(\mu,\nu)$ be the set of all Borel probability measures $\pi$ on $\mbb{R}_+^2$ with marginals $\mu$ and $\nu$, that is,
 \beqnn
\pi(B\times \mbb{R}_+)= \mu(B), ~ \pi(\mbb{R}_+\times B)= \nu(B), \quad B\in \mathscr{B}(\mbb{R_+}).
 \eeqnn
The \textit{$\rho$-Wasserstein distance} $W$ on $\mcr{P}(\mbb{R}_+)$ is defined by
 \beqlb\label{W(mu,nu)def0}
W(\mu,\nu)= \inf_{\pi\in \mcr{C}(\mu,\nu)} \int_{\mbb{R}_+^2} \rho(x,y) \pi(\mathrm{d}x,\mathrm{d}y),
 \quad
\mu,\nu\in \mcr{P}(\mbb{R}_+).
 \eeqlb
It is known that $(\mcr{P}(\mbb{R_+}),W)$ is a complete metric space and the convergence in the distance $W$ is equivalent to the weak convergence of probability measures; see Chen \cite[Theorems~5.4 and~5.6]{Chen04}.

\blemma\label{thW(mu,nu)def1a} Let $\b\mcr{B}_1(\mbb{R_+})$ be the set of functions $f\in \b\mcr{B}(\mbb{R_+})$ satisfying $\|f\|_\infty\le 1$ and $\|f\|_\rho\le 1$. Then we have
 \beqlb\label{W(mu,nu)def1b}
W(\mu,\nu)= \sup_{f\in \b\mcr{B}_1(\mbb{R_+})} |\langle\mu-\nu,f\rangle|,
 \quad
\mu,\nu\in \mcr{P}(\mbb{R}_+).
 \eeqlb
\elemma

\bproof Let $\b\mcr{B}_0(\mbb{R_+})$ be the set of functions $f$ on $\mbb{R}_+$ satisfying $\|f\|_\rho\le 1$. By Chen \cite[Theorem~5.10]{Chen04} it is easy to see that
 \beqnn
W(\mu,\nu)= \sup_{f\in \b\mcr{B}_0(\mbb{R_+})} |\langle\mu-\nu,f\rangle|
 =
\sup_{f\in \b\mcr{B}_0(\mbb{R_+}),f(0)=0} |\langle\mu-\nu,f\rangle|.
 \eeqnn
If $f\in \b\mcr{B}_0(\mbb{R_+})$ satisfies $f(0)=0$, we clearly have
 \beqnn
|f(x)|= |f(x)-f(0)|\le \rho(x,0)\le 1, \quad x\ge 0.
 \eeqnn
and so $f\in \b\mcr{B}_1(\mbb{R_+})$. Then the expression \eqref{W(mu,nu)def1b} follows. \eproof

\blemma\label{thW(mu_1*mu_2,nu_1*nu_2)le.}
For any Borel probability measures $\mu_i$ and $\nu_i$ $(i=1,2)$ on $\mbb{R}_+$, we have
 \beqnn
W(\mu_1*\mu_2,\nu_1*\nu_2)\le W(\mu_1,\nu_1)+W(\mu_2,\nu_2).
 \eeqnn
\elemma

\bproof For $\pi_1\in \mcr{C}(\mu_1,\nu_1)$ and $\pi_2\in \mcr{C}(\mu_2,\nu_2)$, we have $\pi_1*\pi_2\in \mcr{C}(\mu_1*\mu_2,\nu_1*\nu_2)$, and hence
 \beqnn
W(\mu_1*\mu_2,\nu_1*\nu_2)
 \ar\le\ar
\int_{\mbb{R}_+^2}\pi_1(\mathrm{d}x_1,\mathrm{d}y_1)\int_{\mbb{R}_+^2} \big[1\land |(x_1+x_2)-(y_1+y_2)|\big] \pi_2(\mathrm{d}x_2,\mathrm{d}y_2) \cr
 \ar\le\ar
\int_{\mbb{R}_+^2}\pi_1(\mathrm{d}x_1,\mathrm{d}y_1)\int_{\mbb{R}_+^2} (1\land |x_1-y_1| + 1\land |x_2-y_2|) \pi_2(\mathrm{d}x_2,\mathrm{d}y_2) \cr
 \ar=\ar
\int_{\mbb{R}_+^2} (1\land |x_1-y_1|) \pi_1(\mathrm{d}x_1,\mathrm{d}y_1) + \int_{\mbb{R}_+^2} (1\land |x_2-y_2|) \pi_2(\mathrm{d}x_2,\mathrm{d}y_2).
 \eeqnn
Taking the infimum over $\pi_1\in \mcr{C}(\mu_1,\nu_1)$ and $\pi_2\in \mcr{C}(\mu_2,\nu_2)$ gives the desired estimate. \eproof

\subsection{The discrete-time dynamics}\label{SSct2.2}

Let $(\mu_n: n\ge 0)$ be the generalized DR model defined by \eqref{gDR0recur1mu}. Then an corresponding generalized DR process $(X_n: n\ge 0)$ is defined by \eqref{MDR-recurs2X_n}. It is easy to see that, for $n\ge 0$,
 \beqlb\label{MDR-moment1a}
\int_{\mbb{R}_+}z\mu_{n+1}(\d z)\le (1+\alpha m_1)\int_{\mbb{R}_+}z\mu_n(\d z)
 \eeqlb
and
 \beqlb\label{MDR-moment2a}
\int_{\mbb{R}_+}z^2\mu_{n+1}(\d z)
 \ar\le\ar
(1-\alpha)\int_{\mbb{R}_+}z^2\mu_n(\d z) + \alpha\sum_{k=1}^\infty q_k\int_{\mbb{R}_+}z^2\mu_n^{*(k+1)}(\d z) \cr
 \ar\le\ar
(1-\alpha)\int_{\mbb{R}_+}z^2\mu_n(\d z) + \alpha\sum_{k=1}^\infty (k+1)^2q_k\int_{\mbb{R}_+}z^2\mu_n(\d z) \cr
 \ar=\ar
(1+2\alpha m_1 + \alpha m_2) \int_{\mbb{R}_+}z^2\mu_n(\d z).
 \eeqlb
where $m_1$ and $m_2$ denote the first and the second moments of the offspring distribution $q= \{q_1,q_2,\cdots\}$, that is,
 \beqlb\label{m_1,m_2=def}
m_1= \sum_{k=1}^\infty kq_k,
 \quad
m_2= \sum_{k=1}^\infty k^2q_k.
 \eeqlb
For any $f\in \b\mcr{B}(\mbb{R_+})$ we can write
 \beqlb\label{gDRf(X_n)=f(X_0)+sum}
f(X_n)= f(X_0) + \sum_{i=0}^{n-1} A_if(X_i) + M_n(f),
 \eeqlb
where
 \beqlb\label{MDR-A_nf(x)=def1}
A_nf(x)\ar=\ar \alpha\int_{\mbb{R}_+} [f((x+y-1)_+)-f(x)] \mu_n^q(\d y) \cr
 \ar\ar\qqquad
+\, (1-\alpha)[f((x-1)_+)-f(x)] \quad
 \eeqlb
and
 \beqnn
M_n(f)= \sum_{i=0}^{n-1} \big[f(X_{i+1}) - f(X_i) - A_if(X_i)\big].
 \eeqnn
Observe that
 \beqnn
A_if(X_i)\ar=\ar \mbf{E}\big[f\big((x + \eta_iG_i^{-1}(U_i)-1)_+\big) - f(x)\big]\big|_{x=X_i} \ccr
 \ar=\ar
\mbf{E}\big[f\big((X_i + \eta_iG_i^{-1}(U_i)-1)_+\big)\big|X_i\big] - f(X_i) \ccr
 \ar=\ar
\mbf{E}\big[f(X_{i+1}) - f(X_i)\big|X_i\big]
 \eeqnn
and
 \beqnn
M_n(f)= \sum_{i=0}^{n-1} \big\{f(X_{i+1}) - \mbf{E}[f(X_{i+1})|\mcr{F}_i]\big\}.
 \eeqnn
Then $(M_n(f): n\ge 0)$ is a locally bounded martingale.

\subsection{The continuous-time dynamics}\label{SSct2.3}

Recall that $\b\mcr{C}(\mbb{R}_+)$ is the set of bounded continuous functions on $\mbb{R}_+$ and let $\b\mcr{C}^1(\mbb{R}_+)$ be the set of functions in $\b\mcr{C}(\mbb{R}_+)$ with bounded continuous first derivative. Let $\b\mcr{C}_*^1(\mbb{R}_+)$ be the subset of functions $f\in \b\mcr{C}^1(\mbb{R}_+)$ satisfying $f'(0)=0$. For any $f\in \b\mcr{C}_*^1(\mbb{R}_+)$, it is easy to see that
 \beqlb\label{part_tT_tf(x)}
\partial_tT_tf(x)= -f'((x-t)_+)
 =
-T_tf'(x)= -(T_tf)'(x),
 \eeqlb
which is a continuous function of $(t,x)\in \mbb{R}_+^2$. For any Borel function $f$ and any Borel signed-measure $\gamma$ on $\mbb{R}_+$, we write
 \beqnn
\langle\gamma,f\rangle= \int_{\mbb{R}_+} f(x)\gamma(\d x)
 \eeqnn
if the integral exists. Then we may rewrite the differential equation \eqref{Cdifeq0a} more precisely as, for $f\in \b\mcr{C}_*^1(\mbb{R}_+)$,
 \beqlb\label{CDR0difeq1}
\partial_t\langle\mu_t,f\rangle
 =
a\langle\mu_t*\mu_t^q-\mu_t, f\rangle-\langle\mu_t,f'1_{(0,\infty)}\rangle, \quad t\ge 0.
 \eeqlb
Clearly, the above differential equation is equivalent to the integral equation:
 \beqlb\label{CDR(mu)inteq1a}
\langle\mu_t,f\rangle
 =
\langle\mu_0,f\rangle + a\int_0^t\langle\mu_s*\mu_s^q-\mu_s, f\rangle\d s - \int_0^t\langle\mu_s,f'1_{(0,\infty)}\rangle\d s.
 \eeqlb
Moreover, we have the following:

\bproposition\label{CDR(mu)-equiv1a} If the family $(\mu_t: t\ge 0)$ satisfies \eqref{CDR(mu)inteq1a} for every $f\in \b\mcr{C}_*^1(\mbb{R}_+)$, then it satisfies the equation for every $f\in \b\mcr{C}^1(\mbb{R}_+)$. \eproposition

\bproof For each $n\ge 1$ let $r_n\in \b\mcr{C}_*^1(\mbb{R}_+)$ be a function such that $0\le r_n(x)\le n\land x$, $0\le r_n'(x)\le 1$ and $r_n(x)\to l(x):= x$ increasingly as $n\to \infty$ for $x\ge 0$. We can define such a function by
 \beqlb\label{CDR-MPr_n=def}
r_n(x)= \int_0^x g_n(z) \d z, \quad x\ge 0,
 \eeqlb
where
 \beqnn
g_n(z)= \begin{cases}
nz, \ar 0\le z< 1/n,\cr
1, \ar 1/n\le z< n,\cr
n+1-z,~\ar n\le z< n+1,\cr
0, \ar z\ge n+1.
\end{cases}
 \eeqnn
For any $f\in \b\mcr{C}^1(\mbb{R}_+)$, we have $f_n:= f\circ r_n\in \b\mcr{C}_*^1(\mbb{R}_+)$. Then \eqref{CDR(mu)inteq1a} holds for each function $f_n$ by the assumption. By letting $n\to \infty$ and using dominated convergence we see the equation also holds for $f\in \b\mcr{C}^1(\mbb{R}_+)$. \eproof

\bproposition\label{thCequiv1} For a family of probability measures $(\mu_t: t\ge 0)$ on $\mbb{R}_+$, the following properties are equivalent:
 \benumerate

\itm[{\rm(1)}] for every $f\in \b\mcr{C}_*^1(\mbb{R}_+)$ the differential equation \eqref{CDR0difeq1} is satisfied;

\itm[{\rm(2)}] for every $f\in \b\mcr{C}^1(\mbb{R}_+)$ the integral equation \eqref{CDR(mu)inteq1a} is satisfied;

\itm[{\rm(3)}] for every $f\in \b\mcr{B}(\mbb{R}_+)$ the following integral equation is satisfied:
 \beqlb\label{Cinteq1a}
\langle\mu_t,f\rangle= \langle\mu_0,T_tf\rangle + a\int_0^t \langle\mu_s*\mu_s^q-\mu_s, T_{t-s}f\rangle \d s, \quad t\ge 0;
 \eeqlb

\itm[{\rm(4)}] for every $f\in \b\mcr{B}(\mbb{R}_+)$ the following integral equation is satisfied:
 \beqlb\label{Cinteq1b}
\langle\mu_t,f\rangle= \e^{-at}\langle\mu_0,T_tf\rangle + a\int_0^t \e^{a(s-t)} \langle\mu_s*\mu_s^q, T_{t-s}f\rangle \d s, \quad t\ge 0.
 \eeqlb

 \eenumerate
\eproposition

\bproof ``(1)$\Leftrightarrow$(2)'' This follows immediately by Proposition~\ref{CDR(mu)-equiv1a}.

``(1)$\Rightarrow$(3)'' Suppose that $(\mu_t: t\ge 0)$ satisfies \eqref{CDR0difeq1}. For any $f\in \b\mcr{C}_*^1(\mbb{R}_+)$ one can see by \eqref{part_tT_tf(x)} that $(r,t)\mapsto\langle\mu_r, T_tf\rangle$ is continuously differentiable on $[0,\infty)^2$ and, for $t\ge s\ge 0$,
 \beqnn
\frac{\d}{\d s}\langle\mu_s, T_{t-s}f\rangle
 \ar=\ar
\frac{\partial}{\partial r}\langle\mu_r, T_{t-s}f\rangle \Big|_{r=s} - \frac{\partial}{\partial r}\langle\mu_s, T_rf\rangle \Big|_{r=t-s} \cr
 \ar=\ar
\big(a\langle\mu_r*\mu_r^q-\mu_r, T_{t-s}f\rangle - \langle\mu_r,(T_{t-s}f)'\rangle\big)\Big|_{r=s} + \langle\mu_s, T_rf'\rangle \Big|_{r=t-s} \ccr
 \ar=\ar
a\langle\mu_s*\mu_s^q-\mu_s,T_{t-s}f\rangle.
 \eeqnn
Then $(\mu_t: t\ge 0)$ satisfies \eqref{Cinteq1a} for $f\in \b\mcr{C}_*^1(\mbb{R}_+)$. By a monotone class argument we see that \eqref{Cinteq1a} holds for all $f\in \b\mcr{B}(\mbb{R}_+)$.

``(3)$\Rightarrow$(4)'' Suppose that $(\mu_t: t\ge 0)$ satisfies the integral equation \eqref{Cinteq1a}. Then, for any $f\in \b\mcr{C}(\mbb{R}_+)$,
 \beqnn
\frac{\d}{\d s}\langle\mu_s, T_{t-s}f\rangle
 =
a\langle\mu_s*\mu_s^q-\mu_s,T_{t-s}f\rangle,
 \eeqnn
and hence
 \beqnn
\frac{\d}{\d s}\big(\e^{as}\langle\mu_s, T_{t-s}f\rangle\big)
 =
a\e^{as}\langle\mu_s, T_{t-s}f\rangle + \e^s\frac{\d}{\d s}\langle\mu_s, T_{t-s}f\rangle
 =
a\e^{as}\langle\mu_s*\mu_s^q,T_{t-s}f\rangle.
 \eeqnn
By integrating the above equation we get \eqref{Cinteq1b}, which can be extended to $f\in \b\mcr{B}(\mbb{R_+})$.

``(4)$\Rightarrow$(1)'' Suppose that $(\mu_t: t\ge 0)$ satisfies \eqref{Cinteq1b}. For any $f\in \b\mcr{C}_*^1(\mbb{R}_+)$ we can see by \eqref{part_tT_tf(x)} that $t\mapsto \langle\mu_t,f\rangle$ is continuously differentiable and
 \beqnn
\partial_t\langle\mu_t,f\rangle
 \ar=\ar
-\,a\e^{-at}\langle\mu_0,T_tf\rangle - \e^{-at}\langle\mu_0,T_tf'\rangle - a^2\int_0^t \e^{a(s-t)}\langle\mu_s*\mu_s^q-\mu_s, T_{t-s}f\rangle \d s \cr
 \ar\ar
-\, a\int_0^t \e^{a(s-t)}\langle\mu_s*\mu_s^q-\mu_s, T_{t-s}f'\rangle \d s + a\langle\mu_t*\mu_t^q,f\rangle.
 \eeqnn
Using \eqref{Cinteq1b} again we see that $(\mu_t: t\ge 0)$ solves the differential equation \eqref{CDR0difeq1}. \eproof

\bproposition\label{thCsolest1} Suppose that $m_1< \infty$ and $(\mu_t: t\ge 0)$ and $(\gamma_t: t\ge 0)$ are two solutions of \eqref{Cinteq1b}. Then we have
 \beqnn
W(\mu_t,\gamma_t)\le \e^{am_1t} W(\mu_0,\gamma_0), \quad t\ge 0.
 \eeqnn
\eproposition

\bproof Let $f\in \b\mcr{B}_1(\mbb{R_+})$. Then $T_tf\in \b\mcr{B}_1(\mbb{R_+})$ for every $t\ge 0$. Since both $(\mu_t: t\ge 0)$ and $(\gamma_t: t\ge 0)$ are solutions of \eqref{Cinteq1b}, by Lemma~\ref{thW(mu,nu)def1a} we have
 \beqnn
|\langle\mu_t-\gamma_t, f\rangle|
 \ar\le\ar
\e^{-at}|\langle\mu_0-\gamma_0, T_tf\rangle| + a\int_0^t \e^{a(s-t)}|\langle\mu_s*\mu_s^q - \gamma_t*\gamma_s^q, T_{t-s}f\rangle|\d s \cr
 \ar\le\ar
\e^{-at}W(\mu_0,\gamma_0) + a\int_0^t \e^{a(s-t)} W(\mu_s*\mu_s^q,\gamma_t*\gamma_s^q)\d s,
 \eeqnn
where, by Lemma~\ref{thW(mu_1*mu_2,nu_1*nu_2)le.},
 \beqnn
W(\mu_s*\mu_s^q,\gamma_s*\gamma_s^q)
 \le
W(\mu_s,\gamma_s) + W(\mu_s^q,\gamma_s^q)
 \le
(1+m_1)W(\mu_s,\gamma_s).
 \eeqnn
Taking the supremum over all functions $f\in \b\mcr{B}_1(\mbb{R_+})$, we see that
 \beqnn
\e^{at}W(\mu_t,\gamma_t)
 \le
W(\mu_0,\gamma_0) + a(1+m_1)\int_0^t \e^{as} W(\mu_s,\gamma_s)\d s.
 \eeqnn
Then the desired estimate follows by Gronwall's inequality. \eproof.

Now let $\mu_0$ be a fixed probability measure on $\mbb{R}_+$. For $t\ge 0$ define the sub-probability $\mu_t^{(0)}= \e^{-at}\mu_0\circ T_t^{-1}$. Then define the family of sub-probabilities $(\mu_t^{(n)}: t\ge 0)$ for $n\ge 1$ recursively by
 \beqlb\label{CDR-inteq3a}
\langle\mu_t^{(n)},f\rangle= \e^{-at}\langle\mu_0,T_tf\rangle + a\int_0^t \e^{a(s-t)} \langle\mu_s^{(n-1)}*(\mu_s^{(n-1)})^q, T_{t-s}f\rangle \d s.
 \eeqlb

\bproposition\label{thCsolexist1} Suppose that $m_1< \infty$. Then there is a family of probabilities $(\mu_t: t\ge 0)$ on $\mbb{R_+}$ such that
 \beqlb\label{|mu^{(n)}-mu|le..}
\|\mu_t^{(n)}-\mu_t\|_{\var}
 \le
2\sum_{k=n}^\infty\frac{a^{k}(m_1+1)^kt^k}{k!}, \quad t\ge 0,\, n\ge 1.
 \eeqlb
where $\|\cdot\|_{\var} $ denotes the total variation norm. Moreover, the family $(\mu_t: t\ge 0)$ is the unique solution to the integral equation \eqref{Cinteq1b}, where $f\in \b\mcr{B}(\mbb{R_+})$. \eproposition

\bproof The uniqueness of the solution to \eqref{Cinteq1b} holds by Proposition~\ref{thCsolest1}. From \eqref{CDR-inteq3a} it follows that
 \beqnn
|\langle\mu_t^{(n)}-\mu_t^{(n-1)},f\rangle|\ar\le\ar a\int_0^t |\langle\mu_s^{(n-1)}*(\mu_s^{(n-1)})^q - \mu_s^{(n-2)}*(\mu_s^{(n-2)})^q, T_{t-s}f\rangle| \d s \cr
 \ar\le\ar
a\int_0^t \|\mu_s^{(n-1)}*(\mu_s^{(n-1)})^q - \mu_s^{(n-2)}*(\mu_s^{(n-2)})^q\|_{\var} \d s \cr
 \ar\le\ar
a(m_1+1)\int_0^t \|\mu_s^{(n-1)} - \mu_s^{(n-2)}\|_{\var} \d s,
 \eeqnn
where we have used the fact
 \beqnn
\|\mu_1*\nu_1-\mu_2*\nu_2\|_{\var}
 \le
\|\mu_1-\mu_2\|_{\var} + \|\nu_1-\nu_2\|_{\var}.
 \eeqnn
Then for any $0\le t\le u$ we have
 \beqnn
\|\mu_t^{(n)}-\mu_t^{(n-1)}\|_{\var}
 \ar\le\ar
a(m_1+1)\int_0^t \|\mu_{s_1}^{(n-1)}-\mu_{s_1}^{(n-2)}\|_{\var} \d s_1 \cr
 \ar\le\ar
a^2(m_1+1)^2\int_0^t\d s_1\int_0^{s_1} \|\mu_{s_2}^{(n-2)}-\mu_{s_2}^{(n-3)}\|_{\var} \d s_2 \cr
 \ar\le\ar \cdots \cr
 \ar\le\ar
2a^{n-1}(m_1+1)^{n-1}\int_0^t\d s_1 \int_0^{s_1} \cdots \int_0^{s_{n-2}} \d s_{n-1} \ccr
 \ar\le\ar
\frac{2a^{n-1}(m_1+1)^{n-1}t^{n-1}}{(n-1)!},
 \eeqnn
where we have used the fact $\|\mu_{s_{n-1}}^{(1)}-\mu_{s_{n-1}}^{(0)}\|_{\var}\le 2$. Then, for $m> n\ge 1$,
 \beqlb\label{|mu^{(n)}-mu^{(m)}|le..}
\|\mu_t^{(n)}-\mu_t^{(m)}\|_{\var}
 \le
2\sum_{k=n}^{m-1}\frac{a^{k}(m_1+1)^kt^k}{k!}
 \le
2\sum_{k=n}^\infty\frac{a^{k}(m_1+1)^kt^k}{k!}.
 \eeqlb
This shows that $\{\mu_t^{(n)}\}$ is a Cauchy sequence in the total variation distance. Then there are sub-probabilities $(\mu_t: t\ge 0)$ on $\mbb{R_+}$ such that
 \beqnn
\lim_{n\to \infty}\|\mu_t^{(n)}-\mu_t\|_{\var}= 0, \quad t\ge 0.
 \eeqnn
By letting $m\to \infty$ in \eqref{|mu^{(n)}-mu^{(m)}|le..} we obtain \eqref{|mu^{(n)}-mu|le..}. From \eqref{CDR-inteq3a} we see that $(\mu_t: t\ge 0)$ solves \eqref{Cinteq1a} for $f\in \b\mcr{B}(\mbb{R_+})$. In particular, we have
 \beqlb\label{CDR-inteq4a}
\langle\mu_t,1\rangle= \e^{-at}\langle\mu_0,1\rangle + a\int_0^t \e^{a(s-t)} g(\langle\mu_s, 1\rangle) \d s, \quad t\ge 0,
 \eeqlb
where $g$ denotes the probability generating function
 \beqnn
g(z)= \sum_{k=1}^\infty q_kz^{k+1}.
 \eeqnn
Under the assumption $m_1<\infty$, the function $g$ is Lipschitz on $[0,1]$, so $t\mapsto \langle\mu_t,1\rangle\equiv 1$ is the unique solution to the integral equation \eqref{CDR-inteq4a}. Then $(\mu_t: t\ge 0)$ is a family of probabilities. That gives the existence of the solution to \eqref{Cinteq1b}. \eproof

By Propositions~\ref{thCequiv1} and~\ref{thCsolexist1}, the generalized CDR model exists under the condition $m_1< \infty$. By \eqref{Cinteq1b} it is clear that $t\mapsto \langle\mu_t,f\rangle$ is continuous for every $f\in \b\mcr{C}(\mbb{R}_+)$. Then the path $t\mapsto \mu_t$ is continuous by weak convergence of probabilities on $\mbb{R}_+$.

\subsection{A limit theorem for the dynamics}\label{SSct2.4}

Let $a\ge 0$ be a given constant. For each $k\ge a$ let $(\mu_n^{(k)}: n\ge 0)$ be a generalized DR model with renewal rate $\alpha=a/k$ and offspring distribution $q= \{q_1,q_2,\cdots\}$. Let $\gamma_n^{(k)}$ be the probability measure on $\mbb{R}_+$ such that $\gamma_n^{(k)}(\d x)= \mu_n^{(k)}(k\d x)$, $x\ge 0$. By \eqref{gDR0recur1mu} we have
 \beqlb\label{MDR-recurs1a_k}
\gamma_{n+1}^{(k)}= \big[(1-ak^{-1})\gamma_n^{(k)} + ak^{-1}\gamma_n^{(k)} * (\gamma_n^{(k)})^q\big]\circ T_{1/k}^{-1}, \quad n\ge 0.
 \eeqlb

\btheorem\label{thMCDR-est1a} Suppose that $m_1< \infty$. Let $(\mu_t: t\ge 0)$ be the generalized CDR model defined by \eqref{CDR0difeq1}. Then we have
 \beqlb\label{MCDR-est1b}
W(\gamma_{\lfloor kt\rfloor}^{(k)},\mu_t)
 \le
\mathrm{e}^{a(m_1+2)t}\Big[\frac{4}{k}(1+at) + W\big(\gamma^{(k)}_0,\mu_0\big)\Big], \quad t\ge 0.
 \eeqlb
\etheorem

\bproof We first consider an arbitrary function $f\in \b\mcr{B}(\mbb{R_+})$. For any integers $n,n'\ge 0$ satisfying $n+n'= \lfloor kt\rfloor$, we can use \eqref{MDR-recurs1a_k} to see that
 \beqnn
\ar\ar\langle\gamma^{(k)}_{n+1}, T_{(n'-1)/k}f\rangle - \langle\gamma^{(k)}_n, T_{n'/k}f\rangle \ccr
 \ar\ar\qquad
= \big\langle ak^{-1}\gamma_n^{(k)}*(\gamma^{(k)}_n)^q + (1-ak^{-1})\gamma^{(k)}_n, T_{n'/k}f\big\rangle - \langle\gamma^{(k)}_n, T_{n'/k}f\rangle \ccr
 \ar\ar\qquad
=\, ak^{-1}\big\langle\gamma_n^{(k)}*(\gamma^{(k)}_n)^q - \gamma^{(k)}_n, T_{n'/k}f\big\rangle.
 \eeqnn
Summing up the equation over $n$ from $0$ to $\lfloor kt \rfloor-1$ gives
 \beqnn
\langle\gamma^{(k)}_{\lfloor kt \rfloor}, f\rangle
 =
\langle\gamma^{(k)}_0, T_{\lfloor kt \rfloor/k}f\rangle + \frac{a}{k}\sum_{n=0}^{\lfloor kt \rfloor-1}\big\langle\gamma_n^{(k)}*(\gamma^{(k)}_n)^q - \gamma^{(k)}_n, T_{(\lfloor kt \rfloor-n)/k}f\big\rangle.
 \eeqnn
Writing $T^{(k)}_{t,s}= T_{(\lfloor kt \rfloor-\lfloor ks \rfloor)/k}$ for $t\ge s\ge 0$ we obtain
 \beqlb\label{MDR-inteq1a}
\langle\gamma^{(k)}_{\lfloor kt\rfloor}, f\rangle
 =
\langle\gamma^{(k)}_0, T^{(k)}_{t,0}f\rangle + a\int_0^t \big\langle \gamma_{\lfloor ks\rfloor}^{(k)} * (\gamma^{(k)}_{\lfloor ks\rfloor})^q - \gamma^{(k)}_{\lfloor ks\rfloor}, T^{(k)}_{t,s}f\big\rangle\d s + \varepsilon_k(t,f),
 \eeqlb
where
 \beqnn
\varepsilon_k(t,f)= a\Big(t-\frac{\lfloor kt\rfloor}{k}\Big) \big[\langle\gamma^{(k)}_{\lfloor kt \rfloor} - \gamma^{(k)}_{\lfloor kt \rfloor}*(\gamma^{(k)}_{\lfloor kt \rfloor})^q, f\rangle\big].
 \eeqnn
Subtracting \eqref{Cinteq1a} from \eqref{MDR-inteq1a} we get
 \beqlb\label{<gam^(k)_kt,f>-<mu_t,f>}
\langle\gamma^{(k)}_{\lfloor kt\rfloor}, f\rangle - \langle\mu_t, f\rangle
 \ar=\ar
\langle\gamma^{(k)}_0, T^{(k)}_{t,0}f\rangle - \langle\mu_0, T_tf\rangle
-\, a\int_0^t \big(\langle\gamma^{(k)}_{\lfloor ks\rfloor}, T^{(k)}_{t,s}f\rangle - \big\langle\mu_s, T_{t-s}f\big\rangle\big)\d s \cr
 \ar\ar
+\, a\int_0^t \big(\big\langle\gamma^{(k)}_{\lfloor ks\rfloor} * (\gamma^{(k)}_{\lfloor ks\rfloor})^q, T^{(k)}_{t,s}f\big\rangle - \langle\mu_s*\mu_s^q, T_{t-s}f\rangle\big)\d s \ccr
 \ar\ar
+\, \varepsilon_k(t,f).
 \eeqlb
Next we assume $f\in \b\mcr{B}_1(\mbb{R_+})$. Then we also have $T^{(k)}_{t,s}f\in \b\mcr{B}_1(\mbb{R_+})$. For any $r,t\ge 0$ it is easy to see that
 \beqnn
\big|T_tf(x)-T_rf(x)\big|\ar=\ar \big|f((x-t)_+)-f((x-r)_+)\big| \ccr
 \ar\le\ar
\|f\|_\rho\big|(x-t)_+-(x-r)_+\big|\le |t-r|.
 \eeqnn
Then by Lemma~\ref{thW(mu,nu)def1a} we obtain
 \beqnn
\big|\langle\gamma^{(k)}_{\lfloor ks\rfloor}, T^{(k)}_{t,s}f\rangle - \langle\mu_s, T_{t-s}f\rangle\big|
 \ar\le\ar
\langle\gamma^{(k)}_{\lfloor ks\rfloor}, |T^{(k)}_{t,s}f - T_{t-s}f|\rangle + \big|\langle\gamma^{(k)}_{\lfloor ks\rfloor} - \mu_s, T_{t-s}f\rangle\big| \ccr
 \ar\le\ar
\langle\gamma^{(k)}_{\lfloor ks\rfloor}, |T_{(\lfloor kt\rfloor - \lfloor ks\rfloor)/k}f - T_{t-s}f|\rangle + W\big(\gamma^{(k)}_{\lfloor ks\rfloor},\mu_s\big) \ccr
 \ar\le\ar
\frac{2}{k} + W\big(\gamma^{(k)}_{\lfloor ks\rfloor},\mu_s\big).
 \eeqnn
By the same reasoning and an application of lemma \ref{thW(mu_1*mu_2,nu_1*nu_2)le.},
 \beqnn
\big|\big\langle\gamma^{(k)}_{\lfloor ks\rfloor} * (\gamma^{(k)}_{\lfloor ks\rfloor})^q, T^{(k)}_{t,s}f\big\rangle - \langle\mu_s*\mu_s^q, T_{t-s}f\rangle\big|
 \ar\le\ar
\frac{2}{k} + W\big(\gamma^{(k)}_{\lfloor ks\rfloor} * (\gamma^{(k)}_{\lfloor ks\rfloor})^q,\mu_s*\mu_s^q\big) \cr
 \ar\le\ar
\frac{2}{k} + (m_1+1)W\big(\gamma^{(k)}_{\lfloor ks\rfloor},\mu_s\big).
 \eeqnn
For the error term we have
 \beqnn
|\varepsilon_k(t,f)|\le 2\|f\|_\infty\Big(t-\frac{\lfloor kt\rfloor}{k}\Big)\le \frac{2}{k}.
 \eeqnn
With those estimates, by Lemma~\ref{thW(mu,nu)def1a} we deduce from \eqref{<gam^(k)_kt,f>-<mu_t,f>} that
 \beqnn
\big|\langle\gamma^{(k)}_{\lfloor kt\rfloor}, f\rangle - \langle\mu_t, f\rangle\big|
 \le
\frac{4}{k} + W\big(\gamma^{(k)}_0,\mu_0\big) + \frac{4at}{k} + a(m_1+2)\int_0^t W\big(\gamma^{(k)}_{\lfloor ks\rfloor},\mu_s\big) \d s.
 \eeqnn
Taking the supremum over all Lipschitz functions $f\in \b\mcr{B}_1(\mbb{R_+})$ yields
 \beqnn
W\big(\gamma^{(k)}_{\lfloor kt\rfloor},\mu_t\big)
 \le
\frac{4}{k}(1+at) + W\big(\gamma^{(k)}_0,\mu_0\big) + a(m_1+2)\int_0^t W\big(\gamma^{(k)}_{\lfloor ks\rfloor},\mu_s\big) \d s.
 \eeqnn
Then \eqref{MCDR-est1b} follows by Gronwall's inequality. \eproof

As a consequence of Theorem~\ref{thMCDR-est1a} we have the following result:

\btheorem\label{thMDR-lim2a} Suppose that $m_1< \infty$. Let $(\mu_t: t\ge 0)$ be the generalized CDR model defined by \eqref{CDR0difeq1}. If $\gamma_0^{(k)}\overset{\rm w}\to \mu_0$ as $k\to \infty$, then $\gamma_{\lfloor kt\rfloor}^{(k)}\overset{\rm w}\to \mu_t$ for every $t\ge 0$ as $k\to \infty$. \etheorem

\section{The martingale problem}\label{Sect3}

\setcounter{equation}{0}

Let $(\itOmega,\mcr{F},\mbf{P})$ be a complete probability space equipped with a filtration $(\mcr{F}_t: t\ge 0)$ satisfying the usual hypotheses. Let $(A_t: t\ge 0)$ be the family of operators defined by \eqref{CDRA_sf(x)=def}. A positive c\`adl\`ag $(\mcr{F}_t)$-adapted stochastic process is called a solution to the \textit{$(A_t)$-martingale problem} if \eqref{CDR-martprob1b} holds for every $f\in \b\mcr{C}^1(\mbb{R}_+)$, where $\{M_t(f): t\ge 0\}$ is an $(\mcr{F}_t)$-martingale. .

\bproposition\label{CDR-MP3a} If a positive c\`{a}dl\`{a}g $(\mcr{F}_t)$-adapted process $(X_t: t\ge 0)$ satisfies \eqref{CDR-martprob1b} for every $f\in \b\mcr{C}_*^1(\mbb{R}_+)$, then it satisfies \eqref{CDR-martprob1b} for every $f\in \b\mcr{C}^1(\mbb{R}_+)$. \eproposition

\bproof This follows by an approximation of the function $f\in \b\mcr{C}^1(\mbb{R}_+)$ by $f_n:= f\circ r_n\in \b\mcr{C}_*^1(\mbb{R}_+)$, where $r_n$ is given by \eqref{CDR-MPr_n=def}. \eproof

\btheorem\label{CDR-martprob3a} A positive c\`{a}dl\`{a}g process $(X_t: t\ge 0)$ solves the martingale problem \eqref{CDR-martprob1b} if and only if it is a weak solution to the stochastic equation \eqref{gCDRstoeq1X_t}. \etheorem

\bproof If $(X_t: t\ge 0)$ is a weak solution to the stochastic equation \eqref{gCDRstoeq1X_t}, then one can see by It\^o's formula it solves the martingale problem \eqref{CDR-martprob1b}. Conversely, suppose that $(X_t: t\ge 0)$ solves the martingale problem \eqref{CDR-martprob1b}. By It\^o's formula one can see that $Z_t:= \e^{-X_t}$ defines a c\`{a}dl\`{a}g semi-martingale such that
 \beqlb\label{CDR-martprob3b}
Z_t= Z_0 + \int_0^t Z_s1_{\{Z_s<1\}}\d s + a\int_0^t Z_s\d s\int_{\mbb{R}_+} (\e^{-y}-1) \mu_s^q(\d y) + M_t,
 \eeqlb
where $(M_t: t\ge 0)$ is a c\`{a}dl\`{a}g $(\mcr{F}_t)$-martingale. Let $M_1(\d s,\d z)$ be the $(\mcr{F}_t)$-optional time-space random measure on $(0,\infty)\times (\mbb{R}\setminus \{0\})$ defined by
 \beqnn
M_1(\d s,\d z)= \sum_{s>0} 1_{\{\itDelta_s\neq 0\})}\delta_{(s,\itDelta_s)},
 \eeqnn
where $\itDelta_s= M_s-M_{s-}= Z_s-Z_{s-}$. Then we have the orthogonal decomposition
 \beqlb\label{CDR-martprob3c}
M_t= M_0(t) + \int_0^t\int_{\mbb{R}\setminus \{0\}} z\tilde{M}_1(\d s,\d z),
 \eeqlb
where $\{M_0(t): t\ge 0\}$ is a continuous $(\mcr{F}_t)$-martingale and $\tilde{M}_1(\d s,\d z)$ is the compensated measure of $M_1(\d s,\d z)$; see. e.g., \cite[p.353, Theorem~VIII.43]{DeM82}. By \eqref{CDR-martprob3b}, \eqref{CDR-martprob3c} and It\^o's formula, for $f\in \b\mcr{C}^1(\mbb{R}_+)$,
 \beqlb\label{CDR-f(Z_t)=f(Z_0)+(1)}
f(Z_t)\ar=\ar f(Z_0) + a\int_0^t f'(Z_s)Z_s\d s\int_{\mbb{R}_+} (\e^{-y}-1) \mu_s^q(\d y) + \int_0^t f'(Z_s)Z_s1_{\{Z_s<1\}}\d s \cr
 \ar\ar
+ \int_0^t f'(Z_{s-})\d M_0(s) + \int_0^t\int_{\mbb{R}\setminus \{0\}} f'(Z_{s-})y\tilde{M}_1(\d s,\d y) + \frac{1}{2}\int_0^t f''(Z_{s-})\d\langle M_0\rangle(s) \cr
 \ar\ar
+ \int_0^t\int_{\mbb{R}\setminus \{0\}} \big[f(Z_{s-}+y)-f(Z_{s-})-f'(Z_{s-})y\big] M_1(\d s,\d y) \cr
 \ar=\ar
f(Z_0) + a\int_0^t f'(Z_s)Z_s\d s\int_{\mbb{R}_+} (\e^{-y}-1) \mu_s^q(\d y) + \int_0^t f'(Z_s)Z_s1_{\{Z_s<1\}}\d s \cr
 \ar\ar
+\, \text{$(\mcr{F}_t)$-martingale} + \frac{1}{2}\int_0^t f''(Z_{s-})\d\langle M_0\rangle(s) \cr
 \ar\ar
+ \int_0^t\int_{\mbb{R}\setminus \{0\}} \big[f(Z_{s-}+y)-f(Z_{s-})-f'(Z_{s-})y\big] \hat{M}_1(\d s,\d y),
 \eeqlb
where $\hat{M}_1(\d s,\d z)$ is the $(\mcr{F}_t)$-predictable compensator of $M_1(\d s,\d z)$. On the other hand, by applying \eqref{CDRA_sf(x)=def} and \eqref{CDR-martprob1b} directly to the function $x\mapsto f(\e^{-x})$ we see that
 \beqlb\label{CDR-f(Z_t)=f(Z_0)+(2)}
f(Z_t)\ar=\ar f(Z_0) + a\int_0^t \d s\int_{\mbb{R}_+} \big[f(Z_s\e^{-y})-f(Z_s)\big] \mu_s^q(\d y) \cr
 \ar\ar
+ \int_0^t f'(Z_s)Z_s1_{\{Z_s<1\}}\d s + \text{mart.}
 \eeqlb
A comparison of \eqref{CDR-f(Z_t)=f(Z_0)+(1)} and \eqref{CDR-f(Z_t)=f(Z_0)+(2)} shows that
 \beqnn
\ar\ar a\int_0^t\d s\int_{\mbb{R}_+} \big[f(Z_s\e^{-y})-f(Z_s)\big] \mu_s(\d z) \cr
 \ar\ar\qquad
= a\int_0^t f'(Z_s)Z_s\d s\int_{\mbb{R}_+} (\e^{-y}-1) \mu_s(\d y) + \frac{1}{2}\int_0^t f''(Z_s)\d\langle M_0\rangle(s) \cr
 \ar\ar\qquad\quad
+ \int_0^t\int_{\mbb{R}\setminus \{0\}} [f(Z_{s-}+y)-f(Z_{s-})-f'(Z_{s-})y] \hat{M}_1(\d s,\d y).
 \eeqnn
The above equation remains trues for a complex function $f\in \b\mcr{C}^1(\mbb{R}_+)$. In particular, taking $f(x)\equiv \e^{\im\lambda x}$ for $\lambda\in \mbb{R}$ we get
 \beqnn
\ar\ar a\int_0^t\d s\int_{\mbb{R}_+} \big(\e^{\im\lambda Z_s\e^{-y}}-\e^{\im\lambda Z_s}\big) \mu_s^q(\d z) \cr
 \ar\ar\qquad
= \im a\lambda \int_0^t Z_s\e^{\im\lambda Z_s}\d s\int_{\mbb{R}_+} (\e^{-y}-1) \mu_s^q(\d y) - \frac{\lambda^2}{2}\int_0^t \e^{\im\lambda Z_s}\d\langle M_0\rangle(s) \cr
 \ar\ar\qquad\quad
+ \int_0^t\int_{\mbb{R}\setminus \{0\}} \e^{\im\lambda Z_{s-}}(\e^{\im\lambda y}-1-\im\lambda y) \hat{M}_1(\d s,\d y).
 \eeqnn
It follows that
 \beqnn
\ar\ar a\int_0^t\d s\int_{\mbb{R}_+} \e^{\im\lambda Z_s}\big[\e^{\im\lambda Z_s(\e^{-y}-1)} - 1 - \im\lambda Z_s(\e^{-y}-1)\big] \mu_s^a(\d z) \cr
 \ar\ar\qquad
= \int_0^t\int_{\mbb{R}\setminus \{0\}} \e^{\im\lambda Z_{s-}}(\e^{\im\lambda y}-1-\im\lambda y) \hat{M}_1(\d s,\d y) - \frac{\lambda^2}{2}\int_0^t \e^{\im\lambda Z_s}\d\langle M_0\rangle(s),
 \eeqnn
which is an absolutely continuous function of $t\ge 0$. For $T\ge 0$ and $\theta\in \mbb{R}$, integrating the function $t\mapsto \e^{\im(\theta t-\lambda Z_{t-})}$ with respect to both sides over $[0,T]$ we see that
 \beqnn
\ar\ar a\int_0^T\d s\int_{\mbb{R}_+} \e^{\im\theta s}\big[\e^{\im\lambda Z_s(\e^{-y}-1)} - 1 - \im\lambda Z_s(\e^{-y}-1)\big] \mu_s^q(\d y) \cr
 \ar\ar\qquad
= \int_0^T\int_{\mbb{R}\setminus \{0\}} \e^{\im\theta s}(\e^{\im\lambda y}-1-\im\lambda y) \hat{M}_1(\d s,\d y) - \frac{\lambda^2}{2}\int_0^T \e^{\im\theta s}\d\langle M_0\rangle(s).
 \eeqnn
Then the uniqueness of the L\'{e}vy--Khintchine type representation implies that $\langle M_0\rangle(s)\equiv 0$ and, for $t\ge 0$ and $B\in \mcr{B}(\mbb{R}\setminus \{0\})$,
 \beqnn
\hat{M}_1([0,t]\times B)\ar=\ar a\int_0^t\d s\int_{\mbb{R}_+}1_B\big(Z_s(\e^{-y}-1)\big)\mu_s^q(\d y) \cr
 \ar=\ar
a\int_0^t\d s\int_{(0,1)} 1_B\big(Z_{s-}(\e^{-G_s^{-1}(u)}-1)\big)\d u.
 \eeqnn
By a representation theorem, there is a Poisson random measure $N(\d s,\d u)$ on $(0,\infty)\times (0,1)$ with intensity $a\d s\d u$ defined on some extension of the original probability space such that, for $t\ge 0$ and $B\in \mcr{B}(\mbb{R}\setminus \{0\})$,
 \beqnn
M_1([0,t]\times B)= \int_{(0,t]}\int_{(0,1)} 1_B\big(Z_{s-}(\e^{-G_s^{-1}(u)}-1)\big)N(\d s,\d u);
 \eeqnn
see, e.g., \cite[p.93, Theorem 7.4]{IkW89}. Then from \eqref{CDR-martprob3b} it follows that
 \beqnn
Z_t= Z_0 + \int_0^t Z_s1_{\{Z_s<1\}}\d s + \int_{(0,t]}\int_{(0,1)} Z_{s-}(\e^{-G_s^{-1}(u)}-1)N(\d s,\d u).
 \eeqnn
By It\^o's formula one can see that $X_t= -\log Z_t$ is a weak solution of \eqref{gCDRstoeq1X_t}, so it is a DR process associated with $(\mu_t: t\ge 0)$. \eproof

\section{The transition probabilities}\label{Sect4}

\setcounter{equation}{0}

Throughout this section, we fix a generalized CDR model $(\mu_t: t\ge 0)$ defined by \eqref{CDR0difeq1}. For a given constant $r\ge 0$, we are interested in families of probability measures $(\nu_t: t\ge r)$ on $\mbb{R_+}$ solving the differential equation, for $f\in \b\mcr{C}_*^1(\mbb{R_+})$,
 \beqlb\label{Cdifeq2}
\partial_t\langle\nu_t,f\rangle
 =
a\langle\nu_t*\mu_t^q-\nu_t, f\rangle-\langle\nu_t,f'1_{(0,\infty)}\rangle, \quad t\ge r.
 \eeqlb
The above differential equation is equivalent to the integral equation:
 \beqlb\label{CDR(nu)inteq2a}
\langle\nu_t,f\rangle
 =
\langle\nu_r,f\rangle + a\int_r^t \langle\nu_s*\mu_s^q-\nu_s, f\rangle\d s - \int_r^t\langle\nu_s,1_{(0,\infty)}f'\rangle\d s, \quad t\ge r.
 \eeqlb
By arguments similar to those in Subsection~\ref{SSct2.3}, one can prove following results.

\bproposition\label{CDR(nu)-equiv1a} If the family $(\nu_t: t\ge r)$ satisfies \eqref{CDR(nu)inteq2a} for every $f\in \b\mcr{C}_*^1(\mbb{R}_+)$, then it satisfies the equation for every $f\in \b\mcr{C}^1(\mbb{R}_+)$. \eproposition

\bproposition\label{thCequiv2} For a family of probability measures $(\nu_t: t\ge r)$ on $\mbb{R}_+$, the following properties are equivalent:
 \benumerate

\itm[{\rm(1)}] for every $f\in \b\mcr{C}_*^1(\mbb{R}_+)$ the differential equation \eqref{Cdifeq2} is satisfied;

\itm[{\rm(2)}] for every $f\in \b\mcr{C}^1(\mbb{R}_+)$ the integral equation \eqref{CDR(nu)inteq2a} is satisfied;

\itm[{\rm(3)}] for every $f\in \b\mcr{B}(\mbb{R}_+)$ the following integral equation is satisfied:
 \beqlb\label{Cinteq2a}
\langle\nu_t,f\rangle= \langle\nu_r,T_{t-r}f\rangle + a\int_r^t \langle\nu_s*\mu_s^q-\nu_s, T_{t-s}f\rangle \d s, \quad t\ge r;
 \eeqlb

\itm[{\rm(4)}] for every $f\in \b\mcr{B}(\mbb{R}_+)$ the following integral equation is satisfied:
 \beqlb\label{Cinteq2b}
\langle\nu_t,f\rangle= \e^{a(r-t)}\langle\nu_r,T_{t-r}f\rangle + a\int_r^t \e^{a(s-t)} \langle\nu_s*\mu_s^q, T_{t-s}f\rangle \d s, \quad t\ge r.
 \eeqlb

 \eenumerate
\eproposition

\bproposition\label{thCsolest2} Suppose that $(\nu_t: t\ge r)$ and $(\gamma_t: t\ge r)$ are two families of probabilities soling \eqref{Cdifeq2}. Then we have
 \beqlb\label{CDR-solsest1a}
W(\nu_t,\gamma_t)\le \e^{a(t-r)} W(\nu_r,\gamma_r), \quad t\ge r.
 \eeqlb
\eproposition

\bproposition\label{thCsolexist2} To each $\nu_r\in \mcr{P}(\mbb{R}_+)$ there corresponds a unique family of probabilities $(\nu_t: t\ge r)$ solving \eqref{Cinteq2b}. \eproposition

By Proposition~\ref{thCsolexist2}, to each $x\ge 0$ there corresponds a unique family of probabilities $(P_{r,t}(x,\cdot): t\ge r)$ solving the integral equation \eqref{Cinteq2b} with $P_{r,r}(x,\cdot)= \delta_x$. From Proposition~\ref{thCsolest2} it follows that
 \beqlb\label{CDR-W(P_{r,t}(x)P_{r,t}(y))le}
W(P_{r,t}(x,\cdot),P_{r,t}(y,\cdot))\le \e^{a(t-r)} \rho(x,y), \quad t\ge r,\, x,y\in \mbb{R}_+,
 \eeqlb
which implies that the probability measure $P_{r,t}(x,\cdot)$ depends on $x\ge 0$ continuously in the topology of weak convergence. Then $P_{r,t}(x,\cdot)$ is a probability kernel on $\mbb{R}_+$. Given $\gamma\in \mcr{P}(\mbb{R}_+)$, we define $\gamma P_{r,t}\in \mcr{P}(\mbb{R}_+)$ by
 \beqnn
\gamma P_{r,t}(B)= \int_{\mbb{R}_+} \gamma(\d x)P_{r,t}(x,B), \quad B\in \mcr{B}(\mbb{R}_+).
 \eeqnn
It is easy to show that $(\gamma P_{r,t}: t\ge r)$ is the unique solution of \eqref{CDR(nu)inteq2a} with initial state $\gamma$. Consequently, we have
 \beqnn
P_{r,t}(x,\cdot)= \int_{\mbb{R}_+} P_{r,s}(x,\d y)P_{s,t}(x,\cdot), \quad t\ge s\ge r\ge 0.
 \eeqnn
In other words, the family of kernels $(P_{r,t}: t\ge r\ge 0)$ constitute an inhomogeneous Markov transition semigroup on $\mbb{R}_+$. For $f\in \b\mcr{B}(\mbb{R}_+)$, write
 \beqnn
P_{r,t}f(x)= \int_{\mbb{R}_+} f(y)P_{r,t}(x,\d y), \quad t\ge r\ge 0,\,x\in \mbb{R}_+.
 \eeqnn
Then $t\mapsto P_{r,t}f(x)$ is the unique solution to
 \beqlb\label{CDR(P_{r,t})inteq1a}
P_{r,t}f(x)= \e^{a(r-t)}T_{t-r}f(x) + a\int_r^t \e^{a(s-t)}\d s\int_{\mbb{R}_+} \mu_s^q(\d z)\int_{\mbb{R}_+}P_{r,s}(x,\d y)T_{t-s}f(y+z),
 \eeqlb
which is a special case of \eqref{Cinteq2b}. Let $(A_t: t\ge 0)$ be the family of operators defined by \eqref{CDRA_sf(x)=def}. If $f\in \b\mcr{C}^1(\mbb{R}_+)$, then $t\mapsto P_{r,t}f(x)$ is also the unique solution to the \textit{forward integral equation}:
 \beqlb\label{CDR(nu)inteq3a}
P_{r,t}f(x)= f(x) + \int_r^tP_{r,s}A_sf(x)\d s, \quad t\ge r\ge 0,\,x\in \mbb{R}_+,
 \eeqlb
which is a special case of \eqref{CDR(nu)inteq2a}.

\bproposition\label{thCDR(P.f)'=P.f'} For any $t\ge r$ and $f\in \b\mcr{C}_*^1(\mbb{R}_+)$ we have $P_{r,t}f\in \b\mcr{C}^1(\mbb{R}_+)$ and
 \beqlb\label{CDR(P.f)'=P.f'}
(P_{r,t}f)'(x)= P_{r,t}f'(x), \quad t\ge r\ge 0,\,x\in \mbb{R}_+.
 \eeqlb
\eproposition

\bproof For any $f\in \b\mcr{B}(\mbb{R}_+)$ the solution $t\mapsto P_{r,t}f(x)$ to \eqref{CDR(P_{r,t})inteq1a} can be constructed by an iteration argument described as follows. Let $P_{r,t}^{(0)}f(x)= \e^{r-t}T_{t-r}f(x)$. For $n\ge 1$ recursively define
 \beqnn
~~ P_{r,t}^{(n)}f(x)= \e^{a(r-t)}T_{t-r}f(x) + a\int_r^t \e^{a(s-t)}\d s\int_{\mbb{R}_+} \mu_s^q(\d z)\int_{\mbb{R}_+}P_{s,t}^{(n-1)}(x,\d y)T_{t-s}f(y+z).
 \eeqnn
As in the proof of Proposition~\ref{thCsolexist1} one can see that, for $t\ge 0$ and $n\ge 1$,
 \beqnn
\|P_{r,t}^{(n)}f-P_{r,t}f\|_{\infty}
 \le
2\|f\|_{\infty}\sum_{k=n}^\infty\frac{a^{k}(m_1+1)^k(t-r)^k}{k!}.
 \eeqnn
It follows that
 \beqlb\label{lim|P^{(n)}f-Pf|=0}
\lim_{n\to \infty} \|P_{r,t}^{(n)}f-P_{r,t}f\|_{\infty}= 0.
 \eeqlb
For any $f\in \b\mcr{C}_*^1(\mbb{R}_+)$ we have
 \beqnn
(P_{r,t}^{(0)}f)'(x)= \e^{a(r-t)}(T_{t-r}f)'(x)= \e^{a(r-t)}T_{t-r}f'(x)= P_{r,t}^{(0)}f'(x)
 \eeqnn
and, inductively,
 \beqnn
(P_{r,t}^{(n)}f)'(x)= \e^{a(r-t)}T_{t-r}f'(x) + a\int_r^t \e^{a(s-t)}\d s\int_{\mbb{R}_+} \mu_s^q(\d z) \int_{\mbb{R}_+}P_{s,t}^{(n-1)}(x,\d y)T_{t-s}f'(y+z).
 \eeqnn
It follows that $P_{r,t}^{(n)}f\in \b\mcr{C}^1(\mbb{R}_+)$ and $(P_{r,t}^{(n)}f)'= P_{r,t}^{(n)}f'$. By \eqref{lim|P^{(n)}f-Pf|=0} we have
 \beqnn
\lim_{n\to \infty} \|(P_{r,t}^{(n)}f)' - P_{r,t}f'\|_{\infty}
 =
\lim_{n\to \infty} \|P_{r,t}^{(n)}f' - P_{r,t}f'\|_{\infty}= 0,
 \eeqnn
which implies $P_{r,t}f\in \b\mcr{C}^1(\mbb{R}_+)$ and \eqref{CDR(P.f)'=P.f'}. \eproof

\bproposition\label{thCDR(nu)inteq4a} For any $t\ge 0$ and $f\in \b\mcr{C}_*^1(\mbb{R}_+)$ we have the \textit{backward integral equation}:
 \beqlb\label{CDR(nu)inteq4a}
P_{r,t}f(x)= f(x) + \int_r^tA_sP_{s,t}f(x)\d s, \quad 0\le r\le t,\,x\in \mbb{R}_+.
 \eeqlb
\eproposition

\bproof By Proposition~\ref{thCDR(P.f)'=P.f'}, we have $P_{r,t}f\in \b\mcr{C}^1(\mbb{R}_+)$. In view of \eqref{CDR(P_{r,t})inteq1a}, we see that
 \beqlb\label{P_{r-del,r}(x,{0})=}
P_{r,t}(x,\{0\})= P_{r,t}1_{\{0\}}(x)
 =
\e^{r-t}1_{\{x=0\}} + \varepsilon_{r,t}(x).
 \eeqlb
where $\varepsilon_{r,t}(x)\le t-r$. Then for $0< \delta< r$ one can use \eqref{CDR(nu)inteq2a} and the relation $P_{r-\delta,t}f= P_{r-\delta,r}P_{r,t}f$ to see that
 \beqnn
P_{r-\delta,t}f(x)\ar=\ar P_{r,t}f(x) + a\int_{r-\delta}^r\d s\int_{\mbb{R}_+} P_{r-\delta,s}(x,\d y) \int_{\mbb{R}_+} P_{r,t}f(y+z)\mu_s^q(\d z) \cr
 \ar\ar
-\, a\int_{r-\delta}^rP_{r-\delta,s}P_{r,t}f(x)\d s + \int_{r-\delta}^rP_{r-\delta,s}(1_{(0,\infty)}(P_{r,t}f)')(x)\d s \cr
 \ar=\ar
P_{r,t}f(x) + a\int_{r-\delta}^r\d s\int_{\mbb{R}_+} P_{r-\delta,s}(x,\d y) \int_{\mbb{R}_+} P_{r,t}f(y+z)\mu_s^q(\d z) \cr
 \ar\ar
-\, a\int_{r-\delta}^rP_{r-\delta,s}P_{r,t}f(x)\d s + \int_{r-\delta}^rP_{r-\delta,s}(P_{r,t}f)'(x)\d s \cr
 \ar\ar
- \int_{r-\delta}^r [\e^{a(r-t)}1_{\{x=0\}} + \varepsilon_{r-\delta,s}(x)](P_{r,t}f)'(0)\d s.
 \eeqnn
It follows that
 \beqnn
\partial_rP_{r,t}f(x)= aP_{r,t}f(x) - (P_{r,t}f)'(x)1_{\{x>0\}} - a\int_{\mbb{R}_+} P_{r,t}f(x+z)\mu_r^q(\d z)= -A_rP_{r,t}f(x).
 \eeqnn
Then the integral equation \eqref{CDR(nu)inteq4a} holds. \eproof

From \eqref{CDR(nu)inteq3a} and \eqref{CDR(nu)inteq4a} we see that the family of operators $(A_t: t\ge 0)$ is actually a restriction of the \textit{weak generator} of the inhomogeneous transition semigroup $(P_{r,t}: t\ge r\ge 0)$.

\btheorem\label{thCDR-MP5a} A positive c\`{a}dl\`{a}g $(\mcr{F}_t)$-adapted process $(X_t: t\ge 0)$ is a Markov process with inhomogeneous transition semigroup $(P_{r,t}: t\ge r\ge 0)$ if and only if it solves the $(A_t)$-martingale problem. \etheorem

\bproof Suppose that $(X_t: t\ge 0)$ is a Markov process relative to the filtration $(\mcr{F}_t)$ with transition semigroup $(P_{r,t}: t\ge r\ge 0)$. By \eqref{CDR(nu)inteq3a}, for $t\ge r\ge 0$ and $f\in \b\mcr{C}^1(\mbb{R}_+)$ we have
 \beqnn
\mbf{E}[M_t(f)|\mcr{F}_r]
 \ar=\ar
\mbf{E}\Big\{\Big[f(X_t) - f(X_0) - \int_0^t A_sf(X_s) \d s\Big]\Big|\mcr{F}_r\Big\} \cr
 \ar=\ar
\mbf{E}\Big\{\Big[f(X_t) - \int_r^t A_sf(X_s) \d s\Big]\Big|\mcr{F}_r\Big\} - f(X_0) \cr
 \ar\ar
- \int_0^r A_sf(X_s)\d s \cr
 \ar=\ar
P_{r,t}f(X_r) - \int_0^t P_{r,s}A_sf(X_r) \d s - f(X_0) \cr
 \ar\ar
- \int_0^r A_sf(X_s)\d s \cr
 \ar=\ar
f(X_r) - f(X_0) - \int_0^r A_sf(X_s)\d s,
 \eeqnn
which means that $\{M_t(f): t\ge 0\}$ is an $(\mcr{F}_t)$-martingale. Conversely, suppose that for every $f\in \b\mcr{C}^1(\mbb{R}_+)$ the process $\{M_t(f): t\ge 0\}$ defined by \eqref{CDR-martprob1b} is an $(\mcr{F}_t)$-martingale. Then for $v\ge u\ge r\ge 0$ and $F\in \b\mcr{F}_r\subset \b\mcr{F}_u$ we have
 \beqnn
\mbf{E}\{1_F[f(X_v) - f(X_u)]\}
 =
\mbf{E}\bigg\{1_F\int_u^v A_sf(X_s)\d s\bigg\}.
 \eeqnn
Next we assume $f\in \b\mcr{C}_*^1(\mbb{R}_+)$. For $t\ge r\ge 0$, setting $\delta= t-r$, we have
 \beqnn
\ar\ar\mbf{E}\{1_F[f(X_t) - P_{r,t}f(X_r)]\} \ccr
 \ar\ar\qquad
=\mbf{E}\bigg\{1_F\sum_{k=1}^n\big[P_{{r+k\delta/n},t} f(X_{r+k\delta/n}) - P_{{r+(k-1)\delta/n},t} f(X_{r+(k-1)\delta/n})\big]\bigg\} \cr
 \ar\ar\qquad
=\mbf{E}\bigg\{1_F\sum_{k=1}^n\big[P_{{r+k\delta/n},t} f(X_{r+k\delta/n}) - P_{{r+(k-1)\delta/n},t} f(X_{r+k\delta/n})\big]\bigg\} \cr
 \ar\ar\qquad\quad
+\,\mbf{E}\bigg\{1_F\sum_{k=1}^n\big[P_{{r+(k-1)\delta/n},t} f(X_{r+k\delta/n}) - P_{{r+(k-1)\delta/n},t} f(X_{r+(k-1)\delta/n})\big]\bigg\} \cr
 \ar\ar\qquad
= -\,\mbf{E}\bigg\{1_F\sum_{k=1}^n\bigg[\int_{(k-1)\delta/n}^{k\delta/n}
A_{r+s}P_{r+s,t}f(X_{r+k\delta/n})\d s\bigg]\bigg\} \cr
 \ar\ar\qquad\quad
+\,\mbf{E}\bigg\{1_F\sum_{k=1}^n\bigg[\int_{(k-1)\delta/n}^{k\delta/n}
A_{r+s}P_{r+(k-1)\delta/n,t}f(X_{r+s})\d s\bigg]\bigg\} \cr
 \ar\ar\qquad
= -\,\mbf{E}\bigg\{1_F\int_0^\delta A_{r+s}P_{r+s,t}f(X_{r+(\lfloor ns/\delta\rfloor+1)\delta/n}) \d s\bigg\} \cr
 \ar\ar\qquad\quad
+\,\mbf{E}\bigg\{1_F\int_0^\delta A_{r+s}P_{r+\lfloor ns/\delta\rfloor\delta/n,t}f(X_{r+s}) \d s\bigg\}.
 \eeqnn
By \eqref{CDRA_sf(x)=def} and \eqref{CDR(P.f)'=P.f'} one can see that $A_{r+s}P_{r+\lfloor ns/\delta\rfloor\delta/n,t}f(x)\to A_{r+s}P_{r+s,t}(x)$ as $n\to \infty$. Then, by the right continuity of $(X_t: t\ge 0)$, the right-hand side in the above equality tends to zero as $n\to \infty$. It follows that, for $f\in \b\mcr{C}_*^1(\mbb{R}_+)$,
 \beqnn
\mbf{E}\{1_F[f(X_t)]\}= \mbf{E}[1_FP_{r,t}f(X_r)].
 \eeqnn
A monotone class argument shows the above equality holds for all $f\in \b\mcr{B}(\mbb{R}_+)$. This means $(X_t: t\ge 0)$ is a Markov process relative to $(\mcr{F}_t)$ with transition semigroup $(P_{r,t}: t\ge r\ge 0)$. \eproof

By It\^o's formula, one can see that the solution to the stochastic equation \eqref{gCDRstoeq1X_t} also solves the $(A_t)$-martingale problem \eqref{CDR-martprob1b}. Then it is a generalized CDR process by Theorem~\ref{thCDR-MP5a}. Let $\nu_t$ denote the distribution of $X_t$. Then we have $\nu_t= \nu_0P_{0,t}$, so $(\nu_t: t\ge 0)$ is an \textit{entrance law} for the inhomogeneous transition semigroup $(P_{r,t}: t\ge r\ge 0)$. By taking the expectations of the terms in \eqref{CDR-martprob1b}, we see that $(\nu_t: t\ge 0)$ solves \eqref{CDR(nu)inteq2a} for $r=0$. By \eqref{CDR(mu)inteq1a} and the uniqueness of the solution to \eqref{CDR(nu)inteq2a}, if $\nu_0= \mu_0$, then $\nu_t= \mu_t$ for every $t\ge 0$.

\section{A limit theorem for the processes}\label{Sect5}

\setcounter{equation}{0}

In this section, we prove the weak convergence of the rescaled sequence of generalized DR processes in the Skorokhod space. Let $a\ge 0$ be a fixed constant. For each integer $k\ge a$ let $(\mu_n^{(k)}: n\ge 0)$ be a generalized DR model with renewal rate $\alpha=a/k$ and offspring distribution $q= \{q_1,q_2,\cdots\}$, and let $(X_n^{(k)}: n\ge 0)$ be the corresponding generalized DR process. For simplicity, we assume $X_0^{(k)}$ has distribution $\mu_0^{(k)}$. Then $X_n^{(k)}$ has distribution $\mu_n^{(k)}$ for every $n\ge 0$. The process $(X_n^{(k)}: n\ge 0)$ can be constructed recursively by
 \beqlb\label{MDR-recurs^(k)a}
X_{n+1}^{(k)}= \big(X_n^{(k)} + \eta_n^{(k)}(G_n^{(k)})^{-1}(U_n^{(k)}) - 1\big)_+, \quad n\ge 0,
 \eeqlb
where $U_n^{(k)}, \eta_n^{(k)}, X_0^{(k)}$ and $(G_n^{(k)})^{-1}$ are as those in \eqref{MDR-recurs2X_n}, but all depending on the parameter $k$. We shall use the above construction and assume
 \beqlb\label{MDR-mom2cond}
\sup_{k\ge a}\, k^{-2}\mbf{E}\big[(X_0^{(k)})^2\big]
 =
\sup_{k\ge a}\, k^{-2}\int_{\mbb{R}_+} z^2 \mu_0^{(k)}(\d z)< \infty.
 \eeqlb

Let $(A_n^{(k)}: n\ge 0)$ be the generator of $(X_n^{(k)}: n\ge 0)$ and let $(\mcr{F}_n^{(k)}: n\ge 0)$ be its natural filtration. For $k\ge a$ and $f\in \b\mcr{C}^1(\mbb{R}_+)$ write $f_k(x)= f(x/k)$. Then
 \beqlb\label{MDR-seqMP-f_k(X)1}
f_k(X_n^{(k)})= f_k(X_0^{(k)}) + \sum_{i=0}^{n-1} A_i^{(k)}f_k(X_i^{(k)}) + M_n^{(k)}(f_k),
 \eeqlb
where
 \beqnn
A_i^{(k)}f_k(x)\ar=\ar ak^{-1}\int_{\mbb{R}_+} [f_k((x+z-1)_+)-f_k(x)] (\mu_i^{(k)})^q(\d z) \cr
 \ar\ar\qquad
+\, (1-ak^{-1})[f_k((x-1)_+)-f_k(x)]
 \eeqnn
and
 \beqnn
M_n^{(k)}(f_k)\ar=\ar \sum_{i=0}^{n-1} \big\{f_k(X_{i+1}^{(k)}) - \mbf{E}\big[f_k(X_{i+1}^{(k)})\big| \mcr{F}_i^{(k)}\big]\big\}.
 \eeqnn
As observed in Section~2, the process $\{M_n^{(k)}(f_k): n\ge 0\}$ is a locally bounded martingale.

Let $Y_n^{(k)}= X_n^{(k)}/k$ and let $\gamma_n^{(k)}$ be the distribution of $Y_n^{(k)}$. We are interested in the asymptotics of the continuous-time process $(Y_{\lfloor kt\rfloor}^{(k)}: t\ge 0)$ as $k\to \infty$. By \eqref{MDR-seqMP-f_k(X)1} we have
 \beqlb\label{MDR-seqMP-f(Y)1}
f(Y_{\lfloor kt\rfloor}^{(k)})= f(Y_0^{(k)}) + \int_0^{\lfloor kt\rfloor/k} kA_{\lfloor ks\rfloor}^{(k)}f_k(kY_{\lfloor ks\rfloor}^{(k)})\d s + M_{\lfloor kt\rfloor}^{(k)}(f_k),
 \eeqlb
where
 \beqlb\label{MDR-seqkA_[ks]^(k)f_k(ky)1}
kA_{\lfloor ks\rfloor}^{(k)}f_k(ky)\ar=\ar a\int_{\mbb{R}_+} [f((y+z-k^{-1})_+)-f(y)] (\gamma_{\lfloor ks\rfloor}^{(k)})^q(\d z) \cr
 \ar\ar\qquad
+\, (1-ak^{-1})k[f((y-k^{-1})_+)-f(y)]. \quad
 \eeqlb
It is easy to see that
 \beqnn
\big|kA_{\lfloor ks\rfloor}^{(k)}f_k(ky)\big|\le 2a\|f\|_\infty + \|f'\|_\infty.
 \eeqnn
Then $\{M_{\lfloor kt\rfloor}^{(k)}(f_k): t\ge 0\}$ is a locally bounded martingale.

\blemma\label{thMDR-momests} For any $k\ge a$ and $t\ge 0$ we have
 \beqlb\label{MDR-mom12est}
\mbf{E}\big(Y_{\lfloor kt\rfloor}^{(k)}\big)
 \le
\e^{am_1t}\mbf{E}\big(Y_0^{(k)}\big),
 \quad
\mbf{E}\big[\big(Y_{\lfloor kt\rfloor}^{(k)}\big)^2\big]
 \le
\e^{a(2m_2+1)t}\mbf{E}\big[(Y_0^{(k)})^2\big].
 \eeqlb
\elemma

\bproof We only give the proof of the second estimate in \eqref{MDR-mom12est}. The first one follows by similar calculations. For $k\ge a$ and $n\ge 0$ we see from \eqref{MDR-recurs^(k)a} that
 \beqnn
\mbf{E}\big[\big(X_{n+1}^{(k)}\big)^2\big]
 \ar=\ar
\mbf{E}\big[\big(X_n^{(k)}+\eta_n(G_n^{(k)})^{-1}(U_n^{(k)})-1\big)_+^2\big] \ccr
 \ar=\ar
ak^{-1}\int_0^1 \mbf{E}\big[\big(X_n^{(k)}+G_n^{-1}(u)-1\big)_+^2\big]\d u + (1-ak^{-1})\mbf{E}\big[\big(X_n^{(k)}-1\big)_+^2\big] \ccr
 \ar=\ar
ak^{-1}\int_{\mbb{R}_+} \mbf{E}\big[\big(X_n^{(k)}+z-1\big)_+^2\big] (\mu_n^{(k)})^q(\d z) + (1-ak^{-1})\mbf{E}\big[\big(X_n^{(k)}-1\big)_+^2\big] \ccr
 \ar\le\ar
2ak^{-1}\int_{\mbb{R}_+} \big\{\mbf{E}\big[(X_n^{(k)})^2\big]+z^2\big\} (\mu_n^{(k)})^q(\d z) + (1-ak^{-1})\mbf{E}\big[(X_n^{(k)})^2\big] \ccr
 \ar\le\ar
2ak^{-1}\int_{\mbb{R}_+} z^2 (\mu_n^{(k)})^q(\d z) + (1+ak^{-1})\mbf{E}\big[(X_n^{(k)})^2\big] \ccr
 \ar\le\ar
[1+ak^{-1}(2m_2+1)]\mbf{E}\big[(X_n^{(k)})^2\big],
 \eeqnn
where
 \beqnn
\int_{\mbb{R}_+} z^2 (\mu_n^{(k)})^q(\d z)
 \ar=\ar
\sum_{i=1}^\infty q_i\int_{\mbb{R}_+} z^2 (\mu_n^{(k)})^{*i}(\d z) \cr
 \ar=\ar
\sum_{i=1}^\infty i^2q_i\int_{\mbb{R}_+} z^2 \mu_n^{(k)}(\d z)
 =
m_2\mbf{E}\big[(X_n^{(k)})^2\big].
 \eeqnn
It follows that
 \beqnn
\mbf{E}\big[(Y_{\lfloor kt\rfloor}^{(k)})^2\big]
 \le
[1+ak^{-1}(2m_2+1)]^{\lfloor kt\rfloor}\mbf{E}\big[(Y_0^{(k)})^2\big]
 \le
\e^{a(2m_2+1)t}\mbf{E}\big[(Y_0^{(k)})^2\big].
 \eeqnn
That gives the desired estimate. \eproof

\blemma\label{thMDR-mom2sup} For any $k\ge a$ and $t\ge 0$ we have
 \beqlb\label{MDR-mom2sup}
\mbf{E}\Big[\sup_{0\le s\le t}(Y_{\lfloor ks\rfloor}^{(k)})^2\Big]< \infty.
 \eeqlb
\elemma

\bproof Under the assumption \eqref{MDR-mom2cond}, we can use \eqref{MDR-mom12est} to extend \eqref{MDR-seqMP-f(Y)1} to all functions on $\mbb{R}_+$ bounded by $\const\cdot x^2$. In particular, for $k\ge a$ we have
 \beqlb\label{MDR-seqMP-fY1}
Y_{\lfloor kt\rfloor}^{(k)}\ar=\ar Y_0^{(k)} + \int_0^{\lfloor kt\rfloor/k} L_s^{(k)}(Y_{\lfloor ks\rfloor}^{(k)})\d s + M_{\lfloor kt\rfloor}^{(k)},
 \eeqlb
where
 \beqlb\label{MDR-seqL_s^{(k)}(y)=..}
L_s^{(k)}(y)\ar=\ar a\int_{\mbb{R}_+} [(y+z-k^{-1})_+-y] (\gamma_{\lfloor ks\rfloor}^{(k)})^q(\d z) \cr
 \ar\ar\quad
+\, (1-ak^{-1})k[(y-k^{-1})_+-y]
 \eeqlb
and
 \beqlb\label{MDR-seqM_[kt]^(k)=..}
M_{\lfloor kt\rfloor}^{(k)}= \sum_{i=0}^{\lfloor kt\rfloor-1} \big[Y_{i+1}^{(k)} - \mbf{E}\big(Y_{i+1}^{(k)}\big| \mcr{F}_i^{(k)}\big)\big].
 \eeqlb
By \eqref{MDR-seqL_s^{(k)}(y)=..} it is easy to see that
 \beqnn
L_s^{(k)}(y)\le a\int_{\mbb{R}_+} (z-k^{-1})_+ (\gamma_{\lfloor ks\rfloor}^{(k)})^q(\d z)
 \le
a\int_{\mbb{R}_+} z (\gamma_{\lfloor ks\rfloor}^{(k)})^q(\d z)
 =
am_1\mbf{E}\big(Y_{\lfloor ks\rfloor}^{(k)}\big)
 \eeqnn
and
 \beqnn
L_s^{(k)}(y)\ar\ge\ar a\int_{\mbb{R}_+} [(y-k^{-1})_+-y] (\gamma_{\lfloor ks\rfloor}^{(k)})^q(\d z) \cr
 \ar\ar\quad
+\, (1-ak^{-1})k[(y-k^{-1})_+-y]\ge -1.
 \eeqnn
Then, by \eqref{MDR-mom12est},
 \beqlb\label{MDR-seq|L_s(k)(y)|le..}
|L_s^{(k)}(y)|
 \le
1 + am_1\mbf{E}\big(Y_{\lfloor ks\rfloor}^{(k)}\big)
 \le
1 + am_1\e^{am_1s}\mbf{E}(Y_0^{(k)}).
 \eeqlb
Now by \eqref{MDR-seqMP-fY1} and a martingale inequality,
 \beqnn
\mbf{E}\Big[\sup_{0\le s\le t}\big(Y_{\lfloor ks\rfloor}^{(k)}\big)^2\Big]
 \ar\le\ar
3\mbf{E}[(Y_0^{(k)})^2] + 3\mbf{E}\Big[\Big(\int_0^t |L_s^{(k)}(Y_{\lfloor ks\rfloor}^{(k)})|\d s\Big)^2\Big] + 3\mbf{E}\Big[\sup_{0\le s\le t} (M_{\lfloor ks\rfloor}^{(k)})^2\Big] \cr
 \ar\le\ar
3\mbf{E}[(Y_0^{(k)})^2] + 3\Big\{\int_0^t \big[1 + am_1\e^{am_1s}\mbf{E}(Y_0^{(k)})\big] \d s\Big\}^2 \cr
 \ar\ar
+\, 12\mbf{E}\big[(M_{\lfloor kt\rfloor}^{(k)})^2\big].
 \eeqnn
To complete the proof it suffices to show $\mbf{E}[(M_{\lfloor kt\rfloor}^{(k)})^2]< \infty$. By the recursive formula \eqref{MDR-recurs^(k)a}, we have
 \beqnn
Y_{i+1}^{(k)}
 =
(Y_i^{(k)} + k^{-1}\eta_i^{(k)}(G_i^{(k)})^{-1}(U_i^{(k)}) - k^{-1})_+,
 \eeqnn
where $(U_i^{(k)},\eta_i^{(k)})$ is independent of $\mcr{F}_i^{(k)}$. Then
 \beqnn
\mbf{E}(Y_{i+1}^{(k)}\big| \mcr{F}_i^{(k)})
 =
ak^{-1}\int_{\mbb{R}_+} (Y_i^{(k)} + z - k^{-1})_+ (\gamma_i^{(k)})^q(\d z) + (1-ak^{-1})(Y_i^{(k)} - k^{-1})_+.
 \eeqnn
It follows that{\small
 \beqlb\label{MDR-seqE[Y_{i+1}^(k)-E(.)]le..}
\mbf{E}\big\{[Y_{i+1}^{(k)} - \mbf{E}(Y_{i+1}^{(k)}\big| \mcr{F}_i^{(k)})]^2\big\}
 \ar=\ar
ak^{-1} \mbf{E}\Big\{\!\int_{\mbb{R}_+}\!\Big[(Y_i^{(k)} + y - k^{-1})_+\! - (1-ak^{-1})(Y_i^{(k)} - k^{-1})_+ \cr
 \ar\ar\qquad
-\, ak^{-1}\int_{\mbb{R}_+} (Y_i^{(k)} + z - k^{-1})_+ (\gamma_i^{(k)})^q(\d z)\Big]^2 (\gamma_i^{(k)})^q(\d y)\Big\} \cr
 \ar\ar
+\, (1-ak^{-1}) \mbf{E}\Big\{\Big[ak^{-1}(Y_i^{(k)} - k^{-1})_+ \cr
 \ar\ar\qquad
-\, ak^{-1}\int_{\mbb{R}_+} (Y_i^{(k)} + z - k^{-1})_+ (\gamma_i^{(k)})^q(\d z)\Big]^2\Big\} \cr
 \ar\le\ar
ak^{-1}\mbf{E}\Big\{\int_{\mbb{R}_+}\Big[3Y_i^{(k)} + y + \int_{\mbb{R}_+} z (\gamma_i^{(k)})^q(\d z)\Big]^2 (\gamma_i^{(k)})^q(\d y)\Big\} \cr
 \ar\ar\qquad
+\, a^2k^{-2}(1-ak^{-1})\mbf{E}\Big\{\Big[2Y_i^{(k)} + \int_{\mbb{R}_+} z (\gamma_i^{(k)})^q(\d z)\Big]^2\Big\} \cr
 \ar\le\ar
3ak^{-1}\mbf{E}\Big\{\Big[9(Y_i^{(k)})^2 + 2\int_{\mbb{R}_+} z^2 (\gamma_i^{(k)})^q(\d z)\Big]\Big\} \cr
 \ar\ar
+\, 2a^2k^{-2}\mbf{E}\Big\{\Big[4(Y_i^{(k)})^2 + \int_{\mbb{R}_+} z^2 (\gamma_i^{(k)})^q(\d z)\Big]\Big\} \ccr
 \ar\le\ar
ak^{-1}(35+8m_2)\mbf{E}\big[(Y_i^{(k)})^2\big].
 \eeqlb
}By \eqref{MDR-seqM_[kt]^(k)=..}, \eqref{MDR-seqE[Y_{i+1}^(k)-E(.)]le..} and Lemma~\ref{thMDR-momests} we see that
 \beqnn
\mbf{E}\big[(M_{\lfloor kt\rfloor}^{(k)})^2\big]
 \ar=\ar
\sum_{i=0}^{\lfloor kt\rfloor-1} \mbf{E}\big\{\big[Y_{i+1}^{(k)} - \mbf{E}\big(Y_{i+1}^{(k)}\big| \mcr{F}_i^{(k)}\big)\big]^2\big\} \cr
 \ar\le\ar
ak^{-1}(35+8m_2)\sum_{i=0}^{\lfloor kt\rfloor-1} \mbf{E}\big[(Y_i^{(k)})^2\big] \cr
 \ar\le\ar
a(35+8m_2)\int_0^t \mbf{E}\big[(Y_{\lfloor ks\rfloor}^{(k)})^2\big]\d s \cr
 \ar\le\ar
a(35+8m_2)\mbf{E}\big[(Y_0^{(k)})^2\big]\int_0^t \e^{a(2m_2+1)s}\d s< \infty.
 \eeqnn
That proves the desired result. \eproof

\blemma\label{thMDR-incr2est} For $k\ge 1$ let $\tau_k$ be an $(\mcr{F}_{\lfloor kt\rfloor}^{(k)})$-stopping time bounded above by some constant $T\ge 0$. Then for any $t\ge 0$ we have
 \beqnn
\mbf{E}\big[(M_{\lfloor k(\tau_k+t)\rfloor}^{(k)} - M_{\lfloor k\tau_k\rfloor}^{(k)})^2\big]
 \le
a(t+k^{-1})\Big\{35\mbf{E}\Big[\sup_{s\le T+t}(Y_s^{(k)})^2\Big] + 8m_2\e^{a(2m_2+1)(T+t)} \mbf{E}\big[(Y_0^{(k)})^2\big]\Big\}.
 \eeqnn
\elemma

\bproof It is easy to see that both $\lfloor k\tau_k\rfloor$ and $\lfloor k(\tau_k+t)\rfloor$ are stopping times relative to the discrete-time filtration $(\mcr{F}_n^{(k)})$. Write
 \beqnn
M_{\lfloor k(\tau_k+t)\rfloor}^{(k)} - M_{\lfloor k\tau_k\rfloor}^{(k)}
 \ar=\ar
\sum_{i=0}^\infty 1_{\{\lfloor k\tau_k\rfloor+i< \lfloor k(\tau_k+t)\rfloor\}} \big[Y_{\lfloor k\tau_k\rfloor+i+1}^{(k)} - \mbf{E}\big(Y_{\lfloor k\tau_k\rfloor+i+1}^{(k)}\big| \mcr{F}_{\lfloor k\tau_k\rfloor+i}^{(k)}\big)\big].
 \eeqnn
Since $\{\lfloor k\tau_k\rfloor+i< \lfloor k(\tau_k+t)\rfloor\}\in \mcr{F}_{\lfloor k\tau_k\rfloor+i}^{(k)}$, one can show that
 \beqnn
\mbf{E}\big[\big(M_{\lfloor k(\tau_k+t)\rfloor}^{(k)} - M_{\lfloor k\tau_k\rfloor}^{(k)}\big)^2\big]
 \ar=\ar
\sum_{i=0}^\infty \mbf{E}\Big\{1_{\{\lfloor k\tau_k\rfloor+i< \lfloor k(\tau_k+t)\rfloor\}} \big[Y_{\lfloor k\tau_k\rfloor+i+1}^{(k)} \cr
 \ar\ar\qqquad\qquad
-\, \mbf{E}(Y_{\lfloor k\tau_k\rfloor+i+1}^{(k)}\big| \mcr{F}_{\lfloor k\tau_k\rfloor+i}^{(k)})\big]^2\Big\}.
 \eeqnn
By calculations similar to those in \eqref{MDR-seqE[Y_{i+1}^(k)-E(.)]le..} one can see that
 \beqnn
\ar\ar\mbf{E}\Big\{1_{\{\lfloor k\tau_k\rfloor+i< \lfloor k(\tau_k+t)\rfloor\}} \big[Y_{\lfloor k\tau_k\rfloor+i+1}^{(k)} - \mbf{E}(Y_{\lfloor k\tau_k\rfloor+i+1}^{(k)}\big| \mcr{F}_{\lfloor k\tau_k\rfloor+i}^{(k)})\big]^2\Big\} \cr
 \ar\ar\qquad
\le 3ak^{-1}\mbf{E}\Big\{1_{\{\lfloor k\tau_k\rfloor+i< \lfloor k(\tau_k+t)\rfloor\}} \Big[9(Y_{\lfloor k\tau_k\rfloor+i}^{(k)})^2 + 2\int_{\mbb{R}_+} z^2 (\gamma_{\lfloor k\tau_k\rfloor+i}^{(k)})^q(\d z)\Big]\Big\} \cr
 \ar\ar\qquad\quad
+\, 2a^2k^{-2}\mbf{E}\Big\{1_{\{\lfloor k\tau_k\rfloor+i< \lfloor k(\tau_k+t)\rfloor\}} \Big[4(Y_{\lfloor k\tau_k\rfloor+i}^{(k)})^2 + \int_{\mbb{R}_+} z^2 (\gamma_{\lfloor k\tau_k\rfloor+i}^{(k)})^q(\d z)\Big]\Big\} \cr
 \ar\ar\qquad
\le ak^{-1}\mbf{E}\Big\{1_{\{\lfloor k\tau_k\rfloor+i< \lfloor k(\tau_k+t)\rfloor\}} \Big[35\sup_{s\le T+t}(Y_{\lfloor ks\rfloor}^{(k)})^2 + 8\sup_{s\le T+t}\int_{\mbb{R}_+} z^2 (\gamma_{\lfloor ks\rfloor}^{(k)})^q(\d z)\Big]\Big\} \cr
 \ar\ar\qquad
\le ak^{-1}\mbf{E}\Big\{1_{\{\lfloor k\tau_k\rfloor+i< \lfloor k(\tau_k+t)\rfloor\}} \Big[35\sup_{s\le T+t}(Y_{\lfloor ks\rfloor}^{(k)})^2 + 8m_2\e^{a(2m_2+1)(T+t)} \mbf{E}\big[(Y_0^{(k)})^2\big]\Big]\Big\},
 \eeqnn
where we have used Lemma~\ref{thMDR-momests} for the last inequality. It follows that
 \beqnn
\ar\ar\mbf{E}\big[(M_{\lfloor k(\tau_k+t)\rfloor}^{(k)} - M_{\lfloor k\tau_k\rfloor}^{(k)})^2\big] \cr
 \ar\ar\qquad
\le\, ak^{-1}\mbf{E}\Big\{\sum_{i=0}^\infty1_{\{\lfloor k\tau_k\rfloor+i< \lfloor k(\tau_k+t)\rfloor\}} \Big[35\sup_{s\le T+t}(Y_s^{(k)})^2 + 8m_2\e^{a(2m_2+1)(T+t)} \mbf{E}\big[(Y_0^{(k)})^2\big]\Big]\Big\} \cr
 \ar\ar\qquad
=\, ak^{-1}\mbf{E}\Big\{(\lfloor k(\tau_k+t)\rfloor-\lfloor k\tau_k\rfloor) \Big[35\sup_{s\le T+t}(Y_s^{(k)})^2 + 8m_2\e^{a(2m_2+1)(T+t)} \mbf{E}\big[(Y_0^{(k)})^2\big]\Big]\Big\} \cr
 \ar\ar\qquad
\le\, ak^{-1}\mbf{E}\Big\{(kt+1) \Big[35\sup_{s\le T+t}(Y_s^{(k)})^2 + 8m_2\e^{a(2m_2+1)(T+t)} \mbf{E}\big[(Y_0^{(k)})^2\big]\Big]\Big\}.
 \eeqnn
That gives the estimate of the lemma. \eproof

\blemma\label{thMDR-seq0tight} The sequence of processes $\{(Y_{\lfloor kt\rfloor}^{(k)})_{t\ge 0}:$ $k= 1,2,\cdots\}$ is tight in the Skorokhod space $D([0,\infty),\mbb{R}_+)$. \elemma

\bproof For each $k\ge 1$ let $\tau_k$ be an $(\mcr{F}_{\lfloor kt\rfloor}^{(k)})$-stopping time bounded above by some constant $T\ge 0$ and let $\delta_k$ be a constant such that $0\le \delta_k\le 1$ and $\delta_k\to 0$ as $k\to \infty$. From \eqref{MDR-seqMP-fY1} it follows that
 \beqnn
Y_{\lfloor k(\tau_k+\delta_k)\rfloor}^{(k)} - Y_{\lfloor k\tau_k\rfloor}^{(k)}\ar=\ar \int_{\lfloor k\tau_k\rfloor/k}^{\lfloor k(\tau_k+\delta_k)\rfloor/k} L_s^{(k)}(Y_{\lfloor ks\rfloor}^{(k)})\d s + M_{\lfloor k(\tau_k+\delta_k)\rfloor}^{(k)} - M_{\lfloor k\delta_k\rfloor}^{(k)},
 \eeqnn
Then, using \eqref{MDR-seq|L_s(k)(y)|le..},
 \beqnn
\mbf{E}\big[\big(Y_{\lfloor k(\tau_k+\delta_k)\rfloor}^{(k)} - Y_{\lfloor k\delta_k\rfloor}^{(k)}\big)^2\big]
 \ar\le\ar
2\mbf{E}\Big[\Big(\int_{\lfloor k\tau_k\rfloor/k}^{\lfloor k(\tau_k+\delta_k)\rfloor/k} L_s^{(k)}(Y_{\lfloor ks\rfloor}^{(k)})\d s\Big)^2\Big] \cr
 \ar\ar
+\, 2\mbf{E}\Big[\big(M_{\lfloor k(\tau_k+\delta_k)\rfloor}^{(k)} - M_{\lfloor k\delta_k\rfloor}^{(k)}\big)^2\Big] \cr
 \ar\le\ar
2k^{-2}\mbf{E}\Big[\big({\lfloor k(\tau_k+\delta_k)\rfloor} - {\lfloor k\tau_k\rfloor}\big)^2\sup_{0\le s\le T+1} L_s^{(k)}(Y_{\lfloor ks\rfloor}^{(k)})\Big] \cr
 \ar\ar
+\, 2\mbf{E}\Big[\big(M_{\lfloor k(\tau_k+\delta_k)\rfloor}^{(k)} - M_{\lfloor k\delta_k\rfloor}^{(k)}\big)^2\Big] \cr
 \ar\le\ar
2\big(\delta_k+k^{-1}\big)^2\Big[1 + am_1\e^{am_1(T+1)}\mbf{E}(Y_0^{(k)})\Big]^2 \cr
 \ar\ar
+\, 2\mbf{E}\Big[\big(M_{\lfloor k(\tau_k+\delta_k)\rfloor}^{(k)} - M_{\lfloor k\delta_k\rfloor}^{(k)}\big)^2\Big].
 \eeqnn
By Lemma~\ref{thMDR-incr2est}, the right hand side tends to zero as $k\to \infty$. By \eqref{MDR-mom2cond} and Lemma~\ref{thMDR-momests}, the sequence random variables $\{Y_{\lfloor kt\rfloor}^{(k)}:$ $k= 1,2,\cdots\}$ is tight in $\mbb{R}_+$ for each $t\ge 0$. Then the tightness of the sequence of processes $\{(Y_{\lfloor kt\rfloor}^{(k)})_{t\ge 0}:$ $k= 1,2,\cdots\}$ in $D([0,\infty),\mbb{R}_+)$ follows by the result of Aldous \cite[Theorem~1]{Ald78}. \eproof

\btheorem\label{thMDR-conv0Skor} Suppose that $(X_t: t\ge 0)$ is a generalized CDR process associated with the generalized CDR model $(\mu_t: t\ge 0)$, where $X_0$ has distribution $\mu_0$. If the distribution of $Y_0^{(k)}$ converges weakly to $\mu_0$ as $k\to \infty$, then $(Y_{\lfloor kt\rfloor}^{(k)}: t\ge 0)$ converges weakly to $(X_t: t\ge 0)$ in the Skorokhod space $D([0,\infty),\mbb{R}_+)$ as $k\to \infty$. \etheorem

\bproof By Theorem~\ref{thMDR-lim2a} we have $\gamma_{\lfloor kt\rfloor}^{(k)}\overset{\rm w}\to \mu_t$ for every $t\ge 0$ as $k\to \infty$. By Lemma~\ref{thMDR-seq0tight}, the sequence of processes $\{(Y_{\lfloor kt\rfloor}^{(k)})_{t\ge 0}:$ $k= 1,2,\cdots\}$ is tight in the Skorokhod space $D([0,\infty),\mbb{R}_+)$. Then we can pass to a subsequence so that $(Y_{\lfloor kt\rfloor}^{(k)}: t\ge 0)$ converges weakly to some positive c\`adl\`ag process $(X_t: t\ge 0)$ in the topology of $D([0,\infty),\mbb{R}_+)$. By applying the Skorokhod representation, we may assume $(Y_{\lfloor kt\rfloor}^{(k)}: t\ge 0)$ converges a.s.\ to $(X_t: t\ge 0)$ in $D([0,\infty),\mbb{R}_+)$. Clearly, the stochastic equation \eqref{gCDRstoeq1X_t} implies that $\mbf{P}(X= X_{t-})= 1$ for each $t\ge 0$. By \cite[p.118, Proposition~5.2]{EtK86} we have a.s.\ $Y_{\lfloor kt\rfloor}^{(k)}\to X_t$ for each $t\ge 0$. Moreover, since $(X_t: t\ge 0)$ has at most countably many jumps, by \cite[p.118, Proposition~5.2]{EtK86} we also have
 \beqnn
\mbf{P}\big(\text{$Y_{\lfloor kt\rfloor}^{(k)}\to X_t$ for a.e.\ $t\ge 0$}\big)= 1.
 \eeqnn
From \eqref{MDR-seqMP-f(Y)1} it follows that, for $f\in \b\mcr{C}^1_*(\mbb{R}_+)$,
 \beqlb\label{MDR-seqMP-f(Y)2}
f(Y_{\lfloor kt\rfloor}^{(k)})\ar=\ar f(Y_0^{(k)}) + \int_0^{\lfloor kt\rfloor} \big[kA_{\lfloor ks\rfloor}^{(k)}f_k(kY_{\lfloor ks\rfloor}^{(k)})-A_sf(Y_{\lfloor ks\rfloor}^{(k)})\big]\d s \cr
 \ar\ar
+ \int_0^{\lfloor kt\rfloor} A_sf(Y_{\lfloor ks\rfloor}^{(k)})\d s + M_{\lfloor kt\rfloor}^{(k)}(f_k).
 \eeqlb
In view of \eqref{CDRA_sf(x)=def} and \eqref{MDR-seqkA_[ks]^(k)f_k(ky)1}, we can use Lemma~\ref{thW(mu,nu)def1a} and the mean-value theorem to get
 \beqnn
\big|kA_{\lfloor ks\rfloor}^{(k)}f_k(ky)-A_sf(y)\big|
 \ar\le\ar
a\Big|\int_{\mbb{R}_+} [f((y+z-k^{-1})_+)-f(y)] \big[(\gamma_{\lfloor ks\rfloor}^{(k)})^q - \mu_s^q)\big](\d z)\Big| \ccr
 \ar\ar
+\, \big|(1-ak^{-1})k[f((y-k^{-1})_+)-f(y)] + f'(y)\big| \ccr
 \ar\le\ar
am_1\Big|\int_{\mbb{R}_+} [f((y+z-k^{-1})_+)-f(y)] \big(\gamma_{\lfloor ks\rfloor}^{(k)} - \mu_s\big)(\d z)\Big| \ccr
 \ar\ar
+\, \big|k[f((y-k^{-1})_+)-f(y)] + f'(y)\big| \ccr
 \ar\ar
+\, a|f((y-k^{-1})_+)-f(y)| \ccr
 \ar\le\ar
2am_1\big(\|f\|_\infty + \|f'\|_\infty\big) W\big(\gamma_{\lfloor ks\rfloor}^{(k)},\mu_s\big) + \big|f'(y) - f'(\eta)\big| \ccr
 \ar\ar
+\, ak^{-1}\|f'\|_\infty,
 \eeqnn
where $y-k^{-1}\le \eta\le y$. Then, for $f\in \b\mcr{C}^1_*(\mbb{R}_+)$ with uniformly continuous derivative $f'$, by Theorem~\ref{thMCDR-est1a} we have
 \beqnn
\lim_{k\to \infty}\sup_{y\ge 0} \big|kA_{\lfloor ks\rfloor}^{(k)}f_k(ky)-A_sf(y)\big|= 0.
 \eeqnn
From \eqref{MDR-seqMP-f(Y)2} we obtain \eqref{CDR-martprob1b}. By an approximation argument as in the proof of Proposition~\ref{CDR(mu)-equiv1a} we see that \eqref{CDR-martprob1b} holds for $f\in \b\mcr{C}^1(\mbb{R}_+)$. By Theorem~\ref{thCDR-MP5a} we conclude that $(X_t: t\ge 0)$ is a generalized CDR process. Clearly, the arguments above show that any convergent subsequence of $(Y_{\lfloor kt\rfloor}^{(k)}: t\ge 0)$ converges weakly to the generalized CDR process $(X_t: t\ge 0)$ in the space $D([0,\infty),\mbb{R}_+)$ as $k\to \infty$. This gives the desired weak converges. \eproof

\bigskip

\textbf{Acknowledgements.~} This research is supported by the National Key R{\&}D Program of China (No.~2020YFA0712901).

\end{document}